\documentclass[12pt,reqno]{amsart}
\usepackage{amsmath,amssymb,amsfonts,amscd,latexsym,amsthm,mathrsfs,bigints,mathtools,verbatim,cite}
\usepackage[usenames,dvipsnames]{color}
\usepackage[unicode]{hyperref}
\usepackage{hypbmsec}
\newtheorem{theorem}{Theorem}

\theoremstyle{remark}
\newtheorem{remark}{\bf Remark}

\newtheorem*{acknowledgements}{\bf Acknowledgements}

\let\eps\varepsilon

\let\emptyset\varnothing
\renewcommand{\d}{{\mathrm d}}

\renewcommand{\Re}{\operatorname{Re}}

\newcommand{\ord}{\operatorname{ord}}

\newcommand\ba{\boldsymbol a}
\newcommand\bb{\boldsymbol b}

\newcommand\bh{\boldsymbol h}

\newcommand\bx{\boldsymbol x}

\newcommand\bp{\boldsymbol p}
\newcommand\bq{\boldsymbol q}
\newcommand\cF{\mathcal F}
\newcommand\frG{\mathfrak G}
\newcommand\frS{\mathfrak S}

\newcommand\fri{\mathfrak i}
\newcommand\frh{\mathfrak h}
\newcommand\frp{\mathfrak p}
\newcommand\frg{\mathfrak g}

\newcommand\Z{\mathbb Z}

\begin{document}

\title{On cellular rational approximations to $\zeta(5)$}

\author{Francis Brown}
\address{All Souls College, University of Oxford, Oxford, OX1 4AL, UK}
\urladdr{https://www.maths.ox.ac.uk/people/francis.brown}

\author{Wadim Zudilin}
\address{Department of Mathematics, IMAPP, Radboud University, PO Box 9010, 6500~GL Nijmegen, Netherlands}
\urladdr{https://www.math.ru.nl/~wzudilin/}

\date{17 October 2022. \emph{Revised}: 26 January 2026}


%\thanks{}

\begin{abstract}
We analyse a certain family of cellular integrals, which are  period integrals on the moduli space $\mathcal M_{0,8}$ of curves of genus zero with eight marked points, and give  rise to simultaneous rational approximations to $\zeta(3)$ and $\zeta(5)$.
By exploiting the action of a large symmetry group on these integrals, we construct 
%infinitely many \emph{effective} rational approximations $p/q$
an infinite \emph{effective} sequence of rational approximations  $p/q$
to $\zeta(5)$ satisfying
\[
0<\bigg|\zeta(5)-\frac pq\bigg|<\frac1{q^{0.86}}.
\]
\end{abstract}

\maketitle
%==================================================

\section{Introduction}
\label{sec1}

Following Ap\'ery's legendary proof \cite{Ap79} of the irrationality of $\zeta(3)$, the traditional strategy for approaching the irrationality of $\zeta(5)$  is to try to construct linear forms in $1$ and  $\zeta(5)$ with rational coefficients and good arithmetic properties. At the time of writing, it has not succeeded despite many years of effort.  In this paper we propose a  different method for tackling this problem.   It  involves constructing  small linear forms in a larger set of multiple zeta values, which, after setting the unwanted numbers to zero, leads to the approximations described in the abstract. 

The starting point for this method is the recent work \cite{Br16} of one of the authors  which revisits irrationality proofs from a new   geometric perspective. It reproduces  Ap\'ery's irrationality proofs for $\zeta(2)$ and $\zeta(3)$, via Beukers'  famous  re-interpretation \cite{Be79} as Euler-type double and triple integrals,   and produces natural families of multiple integrals of higher order which are linear forms in a  controllable set of multiple zeta values. Whilst a literal generalisation of  Ap\'ery's approach for higher weight zeta values remains a tough challenge, the machinery in \cite{Br16} gives us hope that alternative approaches can be found.

In this manuscript we investigate just one  particular family of such 5-fold integrals, which was  disguised as the family ${}^{\phantom\vee}_8\!\pi^\vee_8$ and highlighted in \cite[Example 7.5]{Br16},
\begin{align}
I(\ba)
&=I(a_1,a_2,a_3,a_4,a_5,a_6,a_7,a_8)
\nonumber\\
&=\idotsint\limits_{0<t_1<t_2<\dots<t_5<1}
\frac{t_1^{a_1}(t_2 - t_1)^{a_2}(t_3 - t_2)^{a_3}(t_4 - t_3)^{a_4}(t_5 - t_4)^{a_5}(1 - t_5)^{a_6}}
{(t_3 - t_1)^{b_{24}}t_3^{b_{14}}(1 - t_4)^{b_{57}}
(t_4 - t_2)^{b_{35}}(t_5 - t_2)^{b_{36}}}
\nonumber\\ &\qquad\qquad\times
\frac{\d t_1\,\d t_2\,\d t_3\,\d t_4\,\d t_5}{(t_3 - t_1)t_3(1 - t_4)(t_4 - t_2)(t_5 - t_2)},
\label{I1}
\end{align}
where
\begin{equation}
\begin{gathered}
b_{24} = a_2 + a_3 + a_6 - a_7 - a_8, \quad b_{14} =  a_4 +  a_7 + a_8  - a_2  - a_6 , \\
b_{57} =  a_4 + a_5 + a_8 - a_2 - a_3, \quad b_{35} = a_2 + a_3 - a_8, \quad b_{36} = a_8
\end{gathered}
\label{b-param}
\end{equation}
and all the parameters $a_1,\dots,a_8$ are assumed to be integers.
According to \cite[Sects.~5.1 and~5.2]{Br16} the integral $I(\ba)$ converges if and only if the following \emph{seventeen} linear forms in the $a_i$ are non-negative:
\begin{equation}
\begin{gathered}
a_1, \; a_2, \; a_3, \; a_4, \; a_5, \; a_6, \; a_7, \; a_1 + a_5 - a_3, \; a_3 + a_6 - a_8, \\
a_4 + a_5 +a_7 + a_8 - a_2 - a_3 - a_6, \; a_7 + a_8 - a_6, \; a_4 + a_8 - a_2, \; a_2 + a_3 + a_6 - a_4 - a_8, \\
a_1 + a_8 - a_3, \; a_1 + a_2 - a_4, \; a_4 + a_5 - a_2, \; a_4 + a_7 +2 a_8      - a_2 - a_3 - a_6.
\end{gathered}
\label{17elem}
\end{equation}
General results on the periods of the  moduli spaces  of genus zero curves with marked points \cite{Br09} imply \emph{a priori} that the family \eqref{I1} of integrals is a $\mathbb Q$-linear combination of multiple zeta values of weight $\leq 5$, namely: 1, $\zeta(2)$, $\zeta(3)$, $\zeta(4)$, $\zeta(5)$ and $\zeta(3) \zeta(2)$.  The cellular nature of this integral  (more precisely, Poincar\'e duality)   suggests that the term of subleading weight  $\zeta(4)$ vanishes, since it is dual to the non-existent `$\zeta(1)$'. Cohomological considerations furthermore imply that the  coefficients  of the two terms  $\zeta(5)$ and $\zeta(2)\zeta(3)$ always occur in the same proportion, i.e.,  there is a single period in leading weight, namely $ \zeta(5) + 2 \zeta(3) \zeta(2)$. 

The main interest of this family is that, additionally, the coefficient of $\zeta(3)$ always vanishes \cite[Sect.~10.2.4]{Br16}. 
Therefore, as  hinted at in \cite{Br16}, the decomposition of $I=I(\ba)$ into a $\mathbb Q$-linear combination of multiple zeta values takes the   very special form:
\begin{equation} 
I=Q\cdot(2\zeta(5)+4\zeta(3)\zeta(2))-4\hat P\cdot\zeta(2)-2P
\label{deco}
\end{equation}
for some $Q,P,\hat P\in\mathbb Q$ (in fact, $Q\in\mathbb Z$ since it may be expressed as a 5-fold residue of the integrand).
The fact that these  linear forms in $\zeta(5)+2\zeta(3)\zeta(2)$, $\zeta(2)$ and 1 are very small  follows from bounds for the integrand along the domain of integration.  From this one may deduce that  at least one of the two numbers $\{\zeta(2), \zeta(5) + 2\zeta(3)\zeta(2)\}$ is irrational, but since this is already known for  $\zeta(2)$, one cannot deduce any new irrationality result from $I$. In this regard, the number $\zeta(2)$ is usually viewed as `parasitic'.

However, another hint from \cite{Br16}  suggests that  the linear forms
\[
I'=I'(\ba)=Q\zeta(5)-P,
\]
obtained by `setting $\zeta(2)$ to zero', are also reasonably small.
 \emph{A priori} this operation does not make sense, but  can be justified  either by cohomological arguments, or by  passing to motivic versions of the integral $I$ and  motivic zeta values, for which it does.

The fact that the $I'$ are small  implies  \emph{a posteriori} that the linear forms
\[
I''=I''(\ba)=Q\zeta(3)-\hat P
\]
are small as well, since $I=2I'+4I''\zeta(2)$. Thus the original cellular integral \eqref{deco} is, in disguise, a pair of simultaneous approximations $I', I''$ to $\zeta(5)$ and $\zeta(3)$. 
It is our principal goal here to quantify these hints from \cite{Br16} as well as to analyse the arithmetic properties of the coefficients $Q=Q(\ba)$, $P=P(\ba)$ and $\hat P=\hat P(\ba)$ of the simultaneous rational approximations  to $\zeta(5)$ and $\zeta(3)$.
As we will see below, there is a large transformation group $\frG$ (of order $7!=5040$) acting on \emph{normalised} versions of the integrals $I(\ba)$ as well as on $I'(\ba),I''(\ba)$ and on the coefficients $Q(\ba),P(\ba),\hat P(\ba)$;
this group allows us to compute  a sharp  denominator  $D=D(\ba) \in \Z$ for which $DI'\in\mathbb Z\zeta(5)+\mathbb Z$.
Such groups famously appear,   and prove themselves to be arithmetically useful, in constructions of rational approximations to $\zeta(2)$, $\zeta(3)$ and $\zeta(4)$ (and other mathematical constants); see \cite{MZ20,RV96,RV01,VZ18,Zu04}.
They are known in the literature by the name of Rhin--Viola groups.
Our `group structure for $\zeta(5)$' shares similarities with its predecessors but also features interesting novelties which we will try to highlight in due course.

One outcome of our construction and analysis is the following result, for which we need to recall a related concept of \emph{effective} rational approximations to a real number $\alpha$, as explained by Nesterenko in his paper~\cite{Ne16}. A family of linear forms $q_n\alpha-p_n$, with $p_n, q_n \in \mathbb{Q}$, is called effective if $p_n,q_n$ are   solutions to a linear difference equation whose coefficients are polynomials in $\mathbb{Q}$  (also known as an Ap\'ery-type recursion). Equivalently, their generating functions satisfy a  differential equation.
We say that a family of rational approximations to $\alpha$ is effective if it can be written in the form $p_n/q_n$, where $p_n,q_n$ are effective. 
Such approximations%
\footnote{Note that passing from $p_n/q_n$ to $p/q$ with $p,q\in\mathbb Z$ may involve multiplying both $p_n$ and $q_n$ by the same rational factor:
rescaling is implicit in the definition.}
$p_n/q_n$ are distinguished from (ineffective!) solutions to
\[
0<\bigg|\alpha-\frac pq\bigg|<\frac1q
\]
in integers $p,q$ whose  existence follows trivially from  squeezing the number in question via the inequalities  $r/q\le\alpha\le(r+1)/q$, with $r=\lfloor q\alpha\rfloor$.

\begin{theorem}
\label{th:z5}
There is an \emph{effective} infinite sequence of  rational approximations $p/q$, with $p,q\in \mathbb{Z}$, to $\zeta(5)$ such that
\[
0<\bigg|\zeta(5)-\frac pq\bigg|<\frac1{q^{0.86}}.
\]
\end{theorem}

Note that this result does not imply the (expected!) irrationality of $\zeta(5)$  which would follow if the upper bound in the inequality were of the form   $1/q^{1+\eps}$ for some $\eps>0$. 
The `worthiness' exponent $0.86$ is nevertheless best possible when compared to any other known constructions of effective rational approximations to $\zeta(5)$
(see \cite{Ne16,Vi22} for a related comparison in the case of Catalan's constant).
In practice, the linear differential equations constructed in this way are  Picard--Fuchs equations  of geometric (Gauss--Manin) origin. 

The result and analysis in this paper also give us confidence that further exploration of the cellular integrals from \cite{Br16} will produce new arithmetic surprises.

A  more general context for this subject is the study of Mellin transforms
\begin{equation*}
%\label{intro: Mellin}
\int_{\gamma}  f_1^{a_1} \dotsb f_n^{a_n} \omega
\end{equation*}
where $f_1,\dots, f_n\colon X \to \mathbb{G}_m$ are morphisms from an algebraic variety  $X$ defined over $\mathbb{Q}$ to the multiplicative  group $\mathbb{G}_m$, the $a_i \in\mathbb{C}$ are complex parameters, $\gamma \subset X(\mathbb{C})$ is a (locally finite)  chain of integration, and $\omega$ is a differential form of degree equal to the dimension of $\gamma$.  In the present note, $X= \mathcal{M}_{0,n}$, the moduli space of Riemann spheres  with $n=8$ marked points, the $f_1,\dots, f_8$ are  cross-ratios, the $a_i$ are integers, and $\omega$ is the 5-form on the second line of~\eqref{I1}.   The combinatorial, arithmetic and analytic structures that we have unearthed in this particular case may  point to the existence of general theorems for  Mellin transforms in algebraic geometry. We hope that our results  may serve as inspiration for future research along these lines.

Since our main task is to emphasise the  ideas and methods behind the analysis of the integrals \eqref{I1}, we have tried to stay non-technical in our exposition as far as possible, and  details of proofs which are either obtainable by finite calculation, or have been verified by computer computation, are omitted. Each section presents different  mathematical features and  structures underlying the integrals \eqref{I1} and ends with some suggestions for generalisation. For the benefit of the reader, 
we have tried to keep our narrative self-contained, and clear and uncluttered by technical details as far as possible. 

\section{Totally symmetric case}
\label{sec2}

We first focus our attention on the `totally symmetric' case when all the parameters $a_1,\dots,a_8$ are equal,
\[
a_1=\dots=a_8=n.
\]
This also means that all the exponents in \eqref{I1} including $b_{24},b_{14},b_{57},b_{35},b_{36}$ are equal to~$n$:
\begin{align*}
I_n=I(n,\dots,n)
&=\idotsint\limits_{0<t_1<t_2<\dots<t_5<1}
\bigg(\frac{t_1(t_2 - t_1)(t_3 - t_2)(t_4 - t_3)(t_5 - t_4)(1 - t_5)}
{(t_3 - t_1)t_3(1 - t_4)(t_4 - t_2)(t_5 - t_2)}\bigg)^n
\\ &\qquad\qquad\times
\frac{\d t_1\,\d t_2\,\d t_3\,\d t_4\,\d t_5}{(t_3 - t_1)t_3(1 - t_4)(t_4 - t_2)(t_5 - t_2)}.
\end{align*}
It is the simplest possible choice of the parameters; as we will witness later, the transformation group $\frG$ acts trivially in this case.

The integrals $I_n$ are effectively computed (up to $n=10$) using Panzer's \texttt{Hyper\-Int}~\cite{Pa15}.
We find out that we indeed have
\begin{equation}
I_n=Q_n\cdot(2\zeta(5)+4\zeta(3)\zeta(2))-4\hat P_n\cdot\zeta(2)-2P_n
\label{I_n}
\end{equation}
for this range; more specifically,
\begin{gather*}
Q_0=1, \; Q_1=21, \; Q_2=2989, \\
\hat P_0=0, \; \hat P_1=\frac{101}4, \; \hat P_2=\frac{344923}{96}, \quad
P_0=0, \; P_1=\frac{87}4, \; P_2=\frac{1190161}{384}
\end{gather*}
for $n=0,1,2$.
Then Koutschan's \texttt{HolonomicFunctions} \cite{Ko10} produces a third order Ap\'ery-type recursion for the integrals $I_n$:
\begin{align*}
&
2(2n + 1)(41218n^3 - 48459n^2 + 20010n - 2871)(n + 1)^5 I_{n+1}
\\ &\;
-(97604224n^9 + 178061760n^8 + 72005308n^7 - 48634688n^6 - 39076836n^5
\\ &\;\quad
+ 2622730n^4 + 7581006n^3 + 920112n^2 - 543402n - 120582) I_n
\\ &\;
-2n(3874492n^8 - 2617900n^7 - 3144314n^6 + 2947148n^5 + 647130n^4 - 1182926n^3
\\ &\;\quad
+ 115771n^2 + 170716n - 44541) I_{n-1}
\\ &\;
+n(41218n^3 + 75195n^2 + 46746n + 9898)(n - 1)^5 I_{n-2}
=0,
\quad\text{where}\; n=2,3,\dots,
\end{align*}
which is also satisfied by the rational coefficients $Q_n,\hat P_n,P_n$.
This already proves the decomposition~\eqref{I_n}, so that $I_n=2I_n'+4I_n''\zeta(2)$ with $I_n'=Q_n\zeta(5)-P_n$ and $I_n''=Q_n\zeta(3)-\hat P_n$.
The characteristic polynomial of the recurrence equation is
\begin{equation*}
4\lambda^3 - 2368\lambda^2 - 188\lambda + 1 \,. 
%\label{charpoly}
\end{equation*}
If
\[
\lambda_1=0.00500378\hdots, \quad \lambda_2=-0.08438431\hdots \quad\text{and}\quad \lambda_3=592.07938053\hdots
\]
denote its roots (ordered according to their absolute value), then a standard localisation procedure leads to the asymptotics
\begin{gather*}
\lim_{n\to\infty}\frac{\log|I_n|}n=\log|\lambda_1|=-5.29756135\dots,
\\
\lim_{n\to\infty}\frac{\log|I_n'|}n=\lim_{n\to\infty}\frac{\log|I_n''|}n=\log|\lambda_2|=-2.47237372\dots
\\ \intertext{and}
\lim_{n\to\infty}\frac{\log|Q_n|}n=\lim_{n\to\infty}\frac{\log|\hat P_n|}n=\lim_{n\to\infty}\frac{\log|P_n|}n
=\log|\lambda_3|=6.38364071\dotsc.
\end{gather*}
Finally, based on an extensive computation of  the rational coefficients $Q_n,\hat P_n,P_n$ we observe experimentally that
\begin{equation}
Q_n, \; d_n^2d_{2n}\hat P_n, \; d_n^5P_n\in\mathbb Z
\quad\text{for}\; n=0,1,2,\dots,
\label{dn-totsym}
\end{equation}
where $d_n$ denotes the least common multiple of $1,2,\dots,n$. As we show later,
\begin{equation}
Q_n=\sum_{k_1=0}^n\binom{n+k_1}{n}{\binom{n}{k_1}}^2
\sum_{k_2=0}^n\binom{n+k_2}{n}{\binom{n}{k_2}}^2
\binom{n+k_1+k_2}{n}
\label{Q_n}
\end{equation}
implying in particular that the coefficients $Q_n$ are integral.
The validity of this formula can be independently established by verifying, again on the basis of \cite{Ko10}, that the double sum on the right-hand side satisfies the above recursion. 

The approximating forms are similar in spirit to the ones for $1,\zeta(2),\zeta(3)$ constructed in \cite[Section~2]{Zu07}
(although  their characteristic polynomials have two roots of  the same size).  A  comment in \emph{loc.\ cit.}\ helps one to identify our approximations with those constructed in \cite[Theorem~1]{Zu02}:
\[
Q_n=\frac{(-1)^{n+1}q_n}{\binom{2n}n}, \;\; \hat P_n=\frac{(-1)^{n+1}\tilde p_n}{\binom{2n}n}, \;\; P_n=\frac{(-1)^{n+1}p_n}{\binom{2n}n}
\quad\text{for}\; n=0,1,2,\dots
\]
(in the notation of \cite{Zu02}). The identities  above may be proven using the  recursions and initial data.

How good are these effective rational approximations in the totally symmetric case? To answer this question, we consider the following measure of 
`worthiness'.  Assume, more generally, that we have constructed \emph{some} sequence of effective rational approximations $q_n\zeta(5)-p_n\in\mathbb Z\zeta(5)+\mathbb Z$ such that
\[
c_0=\lim_{n\to\infty}\frac{\log|q_n\zeta(5)-p_n|}n
\quad\text{and}\quad
c_1=\lim_{n\to\infty}\frac{\log|q_n|}n>c_0.
\]
Then for any choice of $\eps>0$ and all $n$ sufficiently large,
\[
\bigg|\zeta(5)-\frac{p_n}{q_n}\bigg|<\frac1{q_n^{\gamma-\eps}}
\]
with $\gamma=(c_1-c_0)/c_1$.
We call this exponent  $\gamma$ the \emph{worthiness} of the approximations.
When $\gamma>1$, the inequalities imply that $\zeta(5)$ is irrational;
in that case one can also conclude that the irrationality exponent of  $\zeta(5)$ is at most $\gamma/(\gamma-1)$.

In our particular situation here we find that  $c_0=\log|\lambda_2|+5$, and $c_1=\log|\lambda_3|+5$, implying  that the worthiness exponent is $\gamma=0.77795976\dots$: in other words, with the choice $q_n=d_n^5Q_n$, $p_n=d_n^5P_n$ we have
\[
\bigg|\zeta(5)-\frac{p_n}{q_n}\bigg|<\frac1{q_n^{7/9}}
\]
for all $n$ sufficiently large. The group structure for $\zeta(5)$ will enable us to improve this exponent considerably. 

\section{The associated generalised cellular integrals}
\label{sec3}

From \cite{Br16} we know that the subgroup of automorphisms of $\mathcal{M}_{0,8}$ (which is the symmetric group permuting the 8 marked points) preserving the domain of integration in \eqref{I1} is a dihedral group of order~16.
It is generated by a cyclic rotation
\[
\sigma\colon(t_1,t_2,t_3,t_4,t_5)
\mapsto\bigg(1-\frac{t_1}{t_2},1-\frac{t_1}{t_3},1-\frac{t_1}{t_4},1-\frac{t_1}{t_5},1-t_1\bigg)
\]
(of order 8) and a  reflection
\[
\tau\colon(t_1,t_2,t_3,t_4,t_5)
\mapsto\bigg(t_1,\frac{t_1}{t_5},\frac{t_1}{t_4},\frac{t_1}{t_3},\frac{t_1}{t_2}\bigg).
\]
The group however does not preserve the form of the integrand.

Applying to the integral $I(\ba)=I(a_1,a_2,a_3,a_4,a_5,a_6,a_7,a_8)$ the fifth power of~$\sigma$ and passing from  simplicial to cubical coordinates 
\[
t_1=x_1x_2x_3x_4x_5, \;\; t_2=x_2x_3x_4x_5, \;\; t_3=x_3x_4x_5, \;\; t_4=x_4x_5, \;\; t_5=x_5,
\]
we arrive at the integral 
\begin{align}
I(\ba)
&=\mathop{\bigintss\!\!\dotsb\!\!\bigintss}_{[0,1]^5\;\;}
\frac{\splitfrac{x_1^{a_4}(1-x_1)^{a_5}x_2^{a_4 + a_5}(1-x_2)^{a_6}x_3^{a_2 + a_3 + a_6 - a_8}(1-x_3)^{a_1  + a_5 - a_3}}{\times x_4^{a_1 + a_2}(1-x_4)^{a_7}x_5^{a_2}(1-x_5)^{a_1}}}
{\splitfrac{(1-x_1x_2)^{ a_4 + a_5 + a_8 - a_2 - a_3}(1-x_2x_3)^{a_5 + a_6 - a_8}}{\times(1-x_3x_4)^{a_1 + a_2+ a_6 - a_4  - a_8}(1-x_4x_5)^{a_1 + a_7 + a_8 - a_3 - a_6 }}}
\nonumber\\ &\quad\times
\frac{x_2x_3x_4\,\d x_1\,\d x_2\,\d x_3\,\d x_4\,\d x_5}{(1-x_1x_2)(1-x_2x_3)(1-x_3x_4)(1-x_4x_5)}.
\label{I-a}
\end{align}
In this representation, the involution $x_j\mapsto x_{6-j}$ for $j=1,2,3,4,5$ does not change the form of the integral but acts on the set of exponents as follows:
\begin{equation}
\fri_1\colon\ba\mapsto(a_5,a_4,a_3,a_2,a_1,a_7,a_6,  a_4  + a_7 + a_8 - a_2  - a_6  ).
\label{i1a}
\end{equation}
One shows that it generates  a subgroup of order two of the dihedral group alluded to earlier, which preserves the original form of the integrand.

A  substitution
\[
\bx=(x_1,x_2,x_3,x_4,x_5)=
\bigg(\frac{1-y_1}{1-y_1y_2},1-y_1y_2,y_3,1-y_4y_5,\frac{1-y_5}{1-y_4y_5}\bigg)
\]
transforms the 8-parameter integral $I(\ba)$ into a subfamily of  the following 12-parameter  integrals which have  only two rational factors in their denominator:
\begin{align}
J
&=J(\bp;\bq)=J(p_0,p_1,p_2,p_3,p_4,p_5,p_6;q_1,q_2,q_3,q_4,q_5)
\nonumber\\
&=\raisebox{1.8mm}{$\displaystyle\mathop{\bigintsss\!\dotsb\!\bigintsss}_{[0,1]^5\;\;}$}
\frac{\splitfrac{y_1^{p_1}(1-y_1)^{q_1}y_2^{p_2}(1-y_2)^{q_2}y_3^{p_3+1}(1-y_3)^{q_3}y_4^{p_4}(1-y_4)^{q_4}}{\times y_5^{p_5}(1-y_5)^{q_5}\,\d y_1\dotsb\d y_5}}
{(1-y_3(1-y_1y_2))^{p_0+1}(1-y_3(1-y_4y_5))^{p_6+1}}.
\label{J1}
\end{align}
An integral of this form  with 12 parameters reduces to $I(\ba)$ if and only if the   parameters satisfy   the constraints:
\begin{equation}
\begin{gathered}
p_1=p_0+q_4+q_5-q_1-q_3, \;\;
p_3=p_0+q_4+q_5-q_3, \\
p_5=p_0+q_4-q_3, \;\;
p_6=p_0+q_4+q_5-q_1-q_2.
\end{gathered}
\label{cons}
\end{equation}
In this case, the parameters are related  to the $a_1,\ldots, a_8$ as follows: 
\begin{gather*}
p_0=a_5+a_6-a_8, \;\;
p_1=a_2+a_3+a_6-a_4-a_8, \;\;
p_2=a_6, \;\;
p_3=a_2+a_3+a_6-a_8, \\
p_4=a_7, \;\;
p_5=a_3+a_6-a_8, \;\;
p_6=a_1+a_2+a_6-a_4-a_8, \\
q_1=a_4, \;\;
q_2=a_5, \;\;
q_3=a_1+a_5-a_3, \;\;
q_4=a_1, \;\;
q_5=a_2,
\end{gather*}
and, conversely, by 
\begin{gather*}
a_1=q_4, \;\;
a_2=q_5, \;\;
a_3=q_2+q_4-q_3, \;\;
a_4=q_1, \\
a_5=q_2, \;\;
a_6=p_2, \;\;
a_7=p_4, \;\;
a_8=p_2+q_2-p_0.
\end{gather*}

The change of variables $y_j\mapsto y_{6-j}$ for $j=1,2,3,4,5$ shows that the involution \eqref{i1a} naturally lifts to the 12-parameter family as follows:
\begin{equation}
\fri_1\colon(\bp;\bq)\mapsto(p_6,p_5,p_4,p_3,p_2,p_1,p_0;q_5,q_4,q_3,q_2,q_1)
\label{i1pq}
\end{equation}
(the constraints \eqref{cons} are respected by this transformation).
It seems that $\fri_1$ is the only non-trivial automorphism of both families $I(\ba)$ and $J(\bp;\bq)$.

\section{Geometry and cohomology of the cellular integral}
\label{secgeom}
We briefly recall some geometric intuition from \cite{Br16}. A  cellular integral in $\mathcal{M}_{0,8}$ is a certain kind of period of the weight~5 universal mixed Tate motive over $\mathbb{Z}$ defined by   
$H^8( \mathcal{M}_{0,8}^{\delta}, \partial \mathcal{M}_{0,8}^{\delta} )$, and whose weight-graded pieces are direct sums of  $\mathbb{Q}(-n)$,
for $n=0,2,3,4,5$.  Its associated periods are $\mathbb{Q}$-linear combinations of $1,\zeta(2), \zeta(3), \zeta(4), \zeta(5)$ and $\zeta(3)\zeta(2)$. It has a canonical homology class, represented by the domain of integration,  which is invariant under the dihedral group of symmetries of order 16 which acts upon the affine scheme $\mathcal{M}_{0,8}^{\delta}. $

A cellular integral defines a particular subquotient  object $M$ of $ H^8( \mathcal{M}_{0,8}^{\delta}, \partial \mathcal{M}_{0,8}^{\delta} )$.
The integrand of $I(\ba)$ defines a   differential form  and de Rham cohomology class $[\omega_{\ba}] \in M_{dR}$. The divisor of singularities of this form is determined by which out of an  explicit set of  linear forms in the parameters $a_1,\ldots, a_8$ are negative.
%For the  families  that we consider later  on of the form $(a_1,\ldots, a_8) = (\hat a_1n, \ldots, \hat a_8n)$, the motive $M$ only depends on $\hat{\ba}=(\hat a_1,\ldots, \hat a_8)$
%\wadimnote{I change the part to agree with other conventions in the text.}
For the  families $I(\ba n)$ that we consider later, the motive $M$ only depends on $\ba=(a_1,\ldots,a_8)$ but not on~$n$
and can in principle be obtained by a finite, but complicated, computation. We know that $M$  is one-dimensional in highest and lowest weights: $\mathrm{gr}^W_0 M = \mathbb{Q}(0)$, $\mathrm{gr}^W_{10} M = \mathbb{Q}(-5)$, and in the totally symmetric case it is shown in \cite{Br16} that it also vanishes in subleading weight: $\mathrm{gr}^W_8 M =0$.  This calculation can in principle 
be extended to the general non-symmetric case. The dimensions of $M$ in middle weights  have not been rigorously established: to this end, it would be very interesting to have general tools to compute the de Rham realisation $M_{dR}$ via computer. This would  have the considerable benefit of  providing  contiguity relations for the forms $\omega_{\ba}$, and hence direct access to the asymptotics and arithmetic of the forms $I(\ba)$ which we shall deduce later by different methods. 

Computations of periods suggest that the motives in question, for a fixed family
$I(\ba n)$, are of rank~$3$, of the form 
$\operatorname{gr}_W M = \mathbb{Q}(0) \oplus \mathbb{Q}(-2) \oplus \mathbb{Q}(-5)$.   The periods of such a motive are necessarily of the form $1, \zeta(2), \zeta_5=  \mu\zeta(5) + \lambda \zeta(3)\zeta(2)$, for  fixed $\lambda, \mu \in \mathbb{Q}$, which\,---\,in the case of the family of cellular integrals under consideration\,---\,are   $\mu=1, \lambda=2$.
Therefore, a period matrix for $M$ looks like
\begin{equation} \label{periodmatrix1}
\begin{pmatrix} 1 & \zeta(2) & \zeta_5  \\ 
0 & (2\pi i)^2 & (2\pi i)^2 \zeta(3)  \\
0 & 0 & (2\pi i)^5
\end{pmatrix}
\end{equation}
where the columns correspond to a basis  $\omega_0, \omega_2, \omega_5$ of $M_{dR}\cong\mathbb{Q}_{dR}(0) \oplus \mathbb{Q}_{dR}(-2)\oplus \mathbb{Q}_{dR}(-5)$ and the first two rows $\gamma_1, \gamma_2$ correspond to real (i.e., invariant under the action of real Frobenius) homology classes  $M_{B}^{\vee}$, and the last row to an imaginary (real Frobenius anti-invariant) homology class $\gamma_3$.  The integration domain in the integral $I(\ba)$ corresponds to the homology class $\gamma_1$, which is represented by the real simplex $0\leq t_1\leq \dots \leq t_5 \leq 1$. Since the  de Rham class of the integrand  $\omega_{\ba}$ is a $\mathbb{Q}$-linear combination of $\omega_0,\omega_1,\omega_2$, the general integral $I(\ba)$ is a $\mathbb{Q}$-linear combination of the entries in the first row of this period matrix. 
Our linear forms in $1, \zeta(5)$ are obtained by replacing the homology cycle $\gamma_1$ with a particular linear combination of $\gamma_1, \gamma_2$,  in such a way that the  terms  $\zeta(2) = -\frac{1}{24}  (2\pi i)^2 $ cancel out.  This is equivalent to performing a row operation on the first two rows of the matrix \eqref{periodmatrix1}. It necessarily also forces  the term
$\zeta(2)\zeta(3)$  in the top row and right hand column to drop out. This is because the extension group $\operatorname{Ext}^1(\mathbb{Q}(-5), \mathbb{Q}(0))$ in the category of mixed Tate motives over $\mathbb{Z}$ is one-dimensional, with period $\zeta(5)$.

Stated in an equivalent but different way:  there is a change of basis of the Betti homology such that the period matrix takes the form 
\begin{equation} \label{periodmatrix2}
\begin{pmatrix} 1 &  0  & \zeta(5)  \\ 
0 & (2\pi i)^2 & (2\pi i)^2 \zeta(3)  \\
0 & 0 & (2\pi i)^5
\end{pmatrix}.
\end{equation}
This follows from the fact that the  motive $M$ contains at most two non-trivial extensions:  between $\mathbb{Q}(0)$ and $\mathbb{Q}(-5)$, giving a possible period $\zeta(5)$;  and between 
$\mathbb{Q}(-2)$ and $\mathbb{Q}(-5)$, giving a possible period $(2\pi i)^2 \zeta(3)$. Any extension between $\mathbb{Q}(0)$ and $\mathbb{Q}(-2)$ necessarily splits, which explains the shape of the  matrix \eqref{periodmatrix2}, and more precisely  the zero in the top row and middle column. 

In the rest of this note, the word `motivic' means any property of  the integrals $I(\ba)$ which holds on the level of the object $M$, and hence on the entire period matrix. For example, the involution  $\fri_1$  is clearly motivic. It would be very interesting to prove that the entire group  $\frG$ we shall associate  to  $\zeta(5)$  is also motivic.

\section{Barnes-type  representation of the integrals, and asymptotics}
\label{sec4}

The internal integrals over $y_1,y_2$,  and over $y_4,y_5$ in $J(\bp;\bq)$ can be recognised as $_3F_2$-hypergeometric functions using the identity 
\begin{align*}
{}_3F_2\bigg(\begin{matrix} \alpha_0, \, \alpha_1, \, \alpha_2 \\ \beta_1, \, \beta_2 \end{matrix}\biggm|z\bigg)
&=\frac{\Gamma(\beta_1)\Gamma(\beta_2)}{\Gamma(\alpha_1)\Gamma(\alpha_2)\Gamma(\beta_1-\alpha_1)\Gamma(\beta_2-\alpha_2)}
\\ &\quad\times
\iint\limits_{[0,1]^2}\frac{x_1^{\alpha_1-1}(1-x_1)^{\beta_1-\alpha_1-1}x_2^{\alpha_2-1}(1-x_2)^{\beta_2-\alpha_2-1}}
{(1-zx_1x_2)^{\alpha_0}}\,\d x_1\,\d x_2
\end{align*}
and observing that
\[
1-y_3(1-y_1y_2)=\bigg(1+\frac{y_3}{1-y_3}\,y_1y_2\bigg)(1-y_3).
\]
We deduce that
\begin{align}
J(\bp;\bq)
&=\frac{p_1!\,q_1!\,p_2!\,q_2!\,p_4!\,q_4!\,p_5!\,q_5!}{(p_1+q_1+1)!\,(p_2+q_2+1)!\,(p_4+q_4+1)!\,(p_5+q_5+1)!}
\nonumber\\ &\; \qquad  \times
\int_0^1
{}_3F_2\bigg( \begin{matrix} p_0+1, \, p_1+1, \, p_2+1 \\ p_1 + q_1 +2, \, p_2+q_2+2 \end{matrix}\biggm|\frac{-y_3}{1-y_3}\bigg)
\nonumber\\ &\;\quad \qquad \qquad \times
{}_3F_2\bigg( \begin{matrix} p_4+1, \, p_5+1, \, p_6+1 \\ p_4 + q_4 +2, \, p_5+q_5+2 \end{matrix}\biggm|\frac{-y_3}{1-y_3}\bigg)
\,\frac{y_3^{p_3+1}\,\d y_3}{(1-y_3)^{p_0+p_6-q_3+2}}
\nonumber\\ \intertext{which, after changing the variable $z=y_3/(1-y_3)$, becomes }
J(\bp;\bq)
&=\frac{p_1!\,q_1!\,p_2!\,q_2!\,p_4!\,q_4!\,p_5!\,q_5!}{(p_1+q_1+1)!\,(p_2+q_2+1)!\,(p_4+q_4+1)!\,(p_5+q_5+1)!}
\nonumber\\ &\; \qquad \times
\int_0^\infty
{}_3F_2\bigg( \begin{matrix} p_0+1, \, p_1+1, \, p_2+1 \\ p_1 + q_1 +2, \, p_2+q_2+2 \end{matrix}\biggm|-z\bigg)
\nonumber\\ &\;\quad \qquad  \qquad \times
{}_3F_2\bigg( \begin{matrix} p_4+1, \, p_5+1, \, p_6+1 \\ p_4 + q_4 +2, \, p_5+q_5+2 \end{matrix}\biggm|-z\bigg)
\,\frac{z^{p_3+1}\,\d z}{(1+z)^{p_3+q_3-p_0-p_6+1}}.
\label{3F2J}
\end{align}
Applying now the Barnes integral representation
\begin{align*}
&
{}_3F_2\bigg( \begin{matrix} \alpha_0, \, \alpha_1, \, \alpha_2 \\ \beta_1, \, \beta_2 \end{matrix}\biggm|-z\bigg)
\\ &\;
=\frac{\Gamma(\beta_1)\Gamma(\beta_2)}{\Gamma(\alpha_0)\Gamma(\alpha_1)\Gamma(\alpha_2)}
\,\frac1{2\pi i}\int_{-c-i\infty}^{-c+i\infty}\frac{\Gamma(\alpha_1+s)\Gamma(\alpha_2+s)\Gamma(\alpha_3+s)\Gamma(-s)}{\Gamma(\beta_1+s)\Gamma(\beta_2+s)}\,z^s\,\d s,
\end{align*}
where the vertical line $\Re s=-c$ separates the poles of $\Gamma(-s)$ from those of $\Gamma(\alpha_j+s)$ for $j=1,2,3$,
and the Eulerian integral
\[
\int_0^\infty\frac{z^{\alpha-1}}{(1+z)^{\alpha+\beta}}\,\d z
=\frac{\Gamma(\alpha)\Gamma(\beta)}{\Gamma(\alpha+\beta)}
\]
where $\Re\alpha,\Re\beta>0$, we obtain
\begin{align}
J(\bp;\bq)&=\frac{q_1!\,q_2!\,q_4!\,q_5!}{p_0!\,p_6!\,(p_3+q_3-p_0-p_6)!}
\nonumber\\ &\;\times
\frac1{2\pi i}\int_{-c_1-i\infty}^{-c_1+i\infty}\d s
\,\frac{\Gamma(p_0+1+s)\Gamma(p_1+1+s)\Gamma(p_2+1+s)\Gamma(-s)}{\Gamma(p_1+q_1+2+s)\Gamma(p_2+q_2+2+s)}
\nonumber\\ &\;\quad\times
\frac1{2\pi i}\int_{-c_2-i\infty}^{-c_2+i\infty}\d t
\,\frac{\Gamma(p_4+1+t)\Gamma(p_5+1+t)\Gamma(p_6+1+t)\Gamma(-t)}{\Gamma(p_4+q_4+2+t)\Gamma(p_5+q_5+2+t)}
\nonumber\\ &\;\qquad\times
\Gamma(p_3+2+s+t)\Gamma(q_3-p_0-p_6-1-s-t),
\label{intJ}
\end{align}
where the  real numbers $c_1,c_2$ satisfy $0<c_1<1+p_0^*=1+\min\{p_0,p_1,p_2\}$, $0<c_2<1+p_6^*=1+\min\{p_4,p_5,p_6\}$ and $c_1+c_2>1+p_0+p_6-q_3$.

Expanding this Barnes-type double integral into a $\mathbb Q$-linear form in single and double zeta values we find that
\[
J(\bp;\bq)=2Q(\bp;\bq)(\zeta(5)+2\zeta(2)\zeta(3))+\dotsb,
\]
where the extra terms encode a $\mathbb Q$-linear combination of zeta values of weight  strictly less than~5.
Furthermore, the leading coefficient has the following explicit double-sum binomial expression:
\begin{align}
Q(\bp;\bq)
&=(-1)^{p_0+p_1+\dots+p_6}\sum_{k_1,k_2\in\mathbb Z}
\binom{k_1}{p_0}\binom{k_2}{p_6}\binom{k_1+k_2+q_3-p_0-p_6}{p_3+q_3-p_0-p_6}
\nonumber\\ &\quad\times
\binom{q_1}{k_1-p_1}\binom{q_2}{k_1-p_2}\binom{q_4}{k_2-p_4}\binom{q_5}{k_2-p_5}.
\label{sumQ}
\end{align}
In the totally symmetric case (when all $p_j=q_k=n$, except for $p_3=2n$) this gives the explicit expression \eqref{Q_n} displayed in Section~\ref{sec2}.

\begin{remark}
\label{rem-decom}
Performing the decomposition of $J(\bp;\bq)$ into a linear form in multiple zeta values is a difficult technical task.
One potential way of doing so is to find an appropriate collection of contiguouity relations for the integral and perform a related multiple induction in the spirit of \cite{RV96,RV01}.
There is an alternative approach hinted at in \cite{Br09}, which forms the basis of the \texttt{HyperInt} algorithm of \cite{Pa15}: in principle it should be possible to perform the integration  steps over the ring of integers, inverting only those primes which are necessary for taking primitives and performing logarithmic Taylor expansions of hyperlogarithms.  
A strategy we have executed here makes use of the integral \eqref{intJ} and is more in line with \cite{MZ20,Zu04,Zu14}, combined with Beukers' original technology in~\cite{Be79}.
Explicitly, one writes the integrand in \eqref{intJ} as a product of  a rational function and of reciprocals of sines;
decomposes the rational part into a sum of partial fractions and makes relevant shifts of the vertical integration paths. 
The original integral ultimately  becomes a $\mathbb Q$-linear combination of the integrals
\begin{align*}
I_{k_1,k_2}^{(s_1,s_2)}
&=\frac1{(2\pi i)^2}\int_{1/3-i\infty}^{1/3+i\infty}\!\int_{1/3-i\infty}^{1/3+i\infty}
\frac1{(t_1+k_1)^{s_1}(t_2+k_2)^{s_2}}\,
\\ &\qquad\times
\frac\pi{\sin\pi t_1}\,\frac\pi{\sin\pi t_2}\,\frac\pi{\sin\pi(t_1+t_2)}\,\d t_1\,\d t_2
\end{align*}
for $k_1,k_2\in\mathbb Z_{\ge0}$ and $s_1,s_2\in\{1,2\}$.
The very same integrals can then  be  cast in the form
\[
I_{k_1,k_2}^{(s_1,s_2)}
=\frac1{\Gamma(s_1)\Gamma(s_2)}\iint\limits_{[0,1]^2}u^{k_1}v^{k_2}f(u,v)(\log u)^{s_1-1}(\log v)^{s_2-1}\,\frac{\d u}u\,\frac{\d v}v,
\]
where
\[
f(u,v)=\frac1{(2\pi i)^2}\int_{1/3-i\infty}^{1/3+i\infty}\!\int_{1/3-i\infty}^{1/3+i\infty}
u^{t_1}v^{t_2}\,\frac\pi{\sin\pi t_1}\,\frac\pi{\sin\pi t_2}\,\frac\pi{\sin\pi(t_1+t_2)}\,\d t_1\,\d t_2
\]
for $0<u,v<1$. The latter double integral can be explicitly computed via a careful residue analysis:
\[
f(u,v)=\begin{cases}
\dfrac{uv\log u}{(1-u)(1-v)}-\dfrac{u\log(u/v)}{(1-u/v)(1-v)} &\text{if}\; 0<u<v<1,  \\[3pt]
\dfrac{uv\log v}{(1-u)(1-v)}-\dfrac{v\log(v/u)}{(1-u)(1-v/u)} &\text{if}\; 0<v<u<1.
\end{cases}
\]
We skip the  details of the remaining computation and only mention that, for the 8-parameter subfamily subject to~\eqref{cons} translated to the form $I(\ba)$, we indeed find that the latter decomposes as \eqref{deco}, and we have the inclusion \eqref{incl} below  for the  linear form $I'(\ba)$.
A sharp arithmetic analysis  for the companion linear form $I''(\ba)\in\mathbb Z\zeta(3)+\mathbb Q$ is somewhat harder to obtain via the above techniques 
but, fortunately, these linear forms have a different expression which we discuss in Section~\ref{zeta(3)}.
\end{remark}

We now turn our attention to the growth of the quantities associated with $J(\bp;\bq)$. To this end, 
we can adapt standard techniques for computing asymptotics of  Barnes-type integrals (and the coefficients in their decomposition) to the integral in~\eqref{intJ} and expression \eqref{sumQ} (see, e.g., \cite[Sect.~2]{Zu02a} and \cite[Sect.~5]{Zu18}). 
This allows us to compute the asymptotics of $J(\bp n;\bq n)$ as $n\to\infty$ (hence of $I(\ba n)$ as well),
but also of all other integrals and sums that appear in the decomposition of these 5-fold integrals. 
The result of the asymptotic analysis is as follows. Consider a suitable solution $(x,y)$ of the system of equations $F_1=F_2=0$, where
\begin{align*}
%\label{eqn:F1F2}
F_1(x,y)
&=x(p_1+q_1-x)(p_2+q_2-x)(x+y+q_3-p_0-p_6)
\\ &\quad   
-(x-p_0)(x-p_1)(x-p_2)(x+y-p_3),
\\
F_2(x,y)
&=(x+y+q_3-p_0-p_6)(p_4+q_4-y)(p_5+q_5-y)y
\\ &\quad
-(x+y-p_3)(y-p_4)(y-p_5)(y-p_6), 
\end{align*}
then the limit of $|J(\bp n;\bq n)|^{1/n}$ as $n\to\infty$ is equal to
\begin{align*}
&
\frac{|x-p_0|^{p_0}|x-p_1|^{p_1}|x-p_2|^{p_2}|x+y-p_3|^{p_3}}{|p_1+q_1-x|^{p_1+q_1}|p_2+q_2-x|^{p_1+q_2}|x+y+q_3-p_0-p_6|^{p_0+p_6-q_3}}
\\ &\qquad\times
\frac{|y-p_4|^{p_4}|y-p_5|^{p_5}|y-p_6|^{p_6}}{|p_4+q_4-y|^{p_4+q_4}|p_5+q_5-y|^{p_5+q_5}}
\\ &\qquad\times
\frac{q_1^{q_1}q_2^{q_2}q_4^{q_4}q_5^{q_5}}
{p_0^{p_0}(p_3+q_3-p_0-p_6)^{p_3+q_3-p_0-p_6}p_6^{p_6}}.
\end{align*}
If one substitutes other (`non-degenerate') solutions $(x,y)$  into the latter expression, one obtains  asymptotics related to different choices of homology cycles.
This recipe for determining asymptotics may be computed in practice: one can eliminate the variable $y$ by computing the resultant $F(x)$ of the polynomials $F_1$ and $F_2$ with respect to the variable $y$; for any root $x_0$ of $F(x)=0$, the related value of $y$ is found by solving the \emph{linear} equation $F_1(x_0,y)=0$.

When the 12 parameters $p_0,\dots,p_6,q_1,\dots,q_5$ are independent, the resultant $F(x)$ has generic degree~9 in~$x$. When the parameters are subject to
\begin{equation}
p_3+q_3=p_0+q_4+q_5, \;\;
p_3+q_3=p_6+q_1+q_2,
\label{cond12}
\end{equation}
which  are consistent with \eqref{cons},
then the degree of $F(x)$ drops (generically) to~5.
In addition, the system $F_1=F_2=0$ features `trivial' solutions $x=0$, $y=p_3$ and $x=p_3$, $y=0$ when the following two conditions are satisfied:
\[
p_3\in\{p_1+q_1,p_2+q_2\} \quad\text{and}\quad p_3\in\{p_4+q_4,p_5+q_5\}.
\]
These are automatically met under the constraints \eqref{cons} which imply that
\begin{equation}
p_3=p_1+q_1 \quad\text{and}\quad p_3=p_5+q_5.
\label{cons1}
\end{equation}
Thus, under conditions \eqref{cond12}, \eqref{cons1}, which determine when  $J(\bp;\bq)$  reduces to $I(\ba)$, the asymptotics are completely determined by the cubic  polynomial
 \begin{equation} \label{eqn: Fcubic}
\frac{F(x)}{x(p_3-x)}.
\end{equation}
This justifies heuristically the appearance of just \emph{three} (presumably $\mathbb Q$-linearly independent) periods, namely $\zeta(5)+2\zeta(3)\zeta(2)$, $\zeta(2)$ and $1$, in the decomposition of the integrals $I(\ba)$, as well as the rank of the motive~$M$ in Section~\ref{secgeom}. 

\begin{remark}
\label{rem:rec}
The above analysis reveals three \emph{real} asymptotic quantities
\begin{equation}
\lambda_1=\lim_{n\to\infty}I(\ba n)^{1/n},
\quad
\lambda_2=\lim_{n\to\infty}I''(\ba n)^{1/n},
\quad
\lambda_3=\lim_{n\to\infty}Q(\ba n)^{1/n}.
\label{3-asymp}
\end{equation}
The inequality $|\lambda_1|<|\lambda_2|$, which can be checked numerically for any particular admissible choice of $\ba=(a_1,\dots,a_8)$, implies that
\[
\lim_{n\to\infty}I'(\ba n)^{1/n}
=\lim_{n\to\infty}(Q(\ba n)\zeta(5)-P(\ba n))^{1/n}
=\lambda_2,
\]
and implies 
in particular  the nonvanishing of our linear forms $I'(\ba n)\in\mathbb Q+\mathbb Z\zeta(5)$.
A practical (though technically challenging!) task for existing creative telescoping realisations is writing down explicitly a (third order Ap\'ery-type) recursion for the integrals $I(\ba n)=J(\bp n;\bq n)$, which is then automatically valid for their coefficients $Q(\ba n),P(\ba n),\hat P(\ba n)$ and linear forms $I'(\ba n)$.
The numbers \eqref{3-asymp} are precisely the roots of the characteristic polynomial of the recurrence equation (compare with the situation in Section~\ref{sec2}), and the nonvanishing can alternatively be established by showing that the first three entries of $Q(\ba n),P(\ba n),\hat P(\ba n)$ generate three linearly independent solutions (over~$\mathbb C$).
\end{remark}

Finally notice that any automorphism of the 12- or 8-parameter (sub)family $J(\bp;\bq)$ should respect the corresponding asymptotics of $J(\bp n;\bq n)^{1/n}$ as $n\to\infty$. See Section~\ref{sec:Graph} for the results of an analysis along these lines.

\section({Descent to \003\266(3)})%
{Descent to $\zeta(3)$}
\label{zeta(3)}

In this part we assume that conditions \eqref{cond12} and \eqref{cons1} are satisfied, so that we indeed have the decomposition \eqref{deco} for the induced integral $I(\ba)$ and its version $J(\bp;\bq)$.
Notice that formula \eqref{sumQ} is nothing but an iterated residue of the integrand in~\eqref{J1}.
Our computation mentioned in Remark~\ref{rem-decom} reveals that $I''(\ba)$ can be read off from
\[
\frac1{(2\pi i)^2}\raisebox{1.8mm}{$\displaystyle\mathop{\bigintsss\!\dotsb\!\bigintsss}_{[0,1]^3,\;|y_4|=|y_5|=\eps}$}
\frac{\splitfrac{y_1^{p_1}(1-y_1)^{q_1}y_2^{p_2}(1-y_2)^{q_2}y_3^{p_3+1}(1-y_3)^{q_3}y_4^{p_4}(1-y_4)^{q_4}}{\times y_5^{p_5}(1-y_5)^{q_5}\,\d y_1\dotsb\d y_5}}
{(1-y_3(1-y_1y_2))^{p_0+1}(1-y_3(1-y_4y_5))^{p_6+1}}
\]
for $\eps<1/4$; one can also choose $|y_1|=|y_2|=\eps$ instead. Using
\begin{align*}
\frac1{(1-y_3(1-y_4y_5))^{p_6+1}}
&=\frac1{(y_3y_4y_5)^{p_6+1}}\bigg(1+\frac{1-y_3}{y_3y_4y_5}\bigg)^{-(p_6+1)}
\\
&=\frac1{(y_3y_4y_5)^{p_6+1}}\sum_{k\ge p_6}\binom{k}{p_6}\bigg(-\frac{1-y_3}{y_3y_4y_5}\bigg)^{k-p_6}
\\
&=\sum_{k\in\mathbb Z}(-1)^{k+p_6}\binom{k}{p_6}(y_3y_4y_5)^{-k-1}(1-y_3)^{-p_6+k},
\end{align*}
so that
\begin{align*}
&
\frac{y_3^{p_3+1}(1-y_3)^{q_3}y_4^{p_4}(1-y_4)^{q_4}y_5^{p_5}(1-y_5)^{q_5}}
{(1-y_3(1-y_4y_5))^{p_6+1}}
\\ &\qquad
=(-1)^{p_6}\sum_{k\in\mathbb Z}(-1)^k\binom{k}{p_6}y_3^{p_3-k}(1-y_3)^{q_3-p_6+k}
\frac{(1-y_4)^{q_4}}{y_4^{k-p_4+1}}\,\frac{(1-y_5)^{q_5}}{y_5^{k-p_5+1}},
\end{align*}
and applying the following formula with $y=y_4$ and $y=y_5$:
\[
\frac1{2\pi i}\oint_{|y|=\eps}\frac{(1-y)^q}{y^{k-p+1}}\,\d y
=\frac1{2\pi i}\oint_{|y|=\eps}\sum_{m=0}^q\binom qm(-1)^m\frac{y^m\,\d y}{y^{k-p+1}}
=(-1)^{k-p}\binom q{k-p},
\]
we find out that
\begin{align}
I''(\ba)
&=(-1)^{p_4+p_5+p_6}\sum_{k\in\mathbb Z}(-1)^k\binom{k}{p_6}\binom{q_4}{k-p_4}\binom{q_5}{k-p_5}
\nonumber\\ &\qquad\times
\iiint\limits_{[0,1]^3}
\frac{y_1^{p_1}(1-y_1)^{q_1}y_2^{p_2}(1-y_2)^{q_2}y_3^{p_3-k}(1-y_3)^{q_3-p_6+k}\,\d y_1\,\d y_2\,\d y_3}
{(1-y_3(1-y_1y_2))^{p_0+1}}
\nonumber\displaybreak[2]\\
&=(-1)^{p_4+p_5+p_6}\sum_{k\in\mathbb Z}(-1)^k\binom{k}{p_6}\binom{q_4}{k-p_4}\binom{q_5}{k-p_5}
\nonumber\\ &\qquad\times
J_3(p_0,p_1,p_2,p_3-k;q_1,q_2,q_3-p_6+k),
\label{I3}
\end{align}
where
\begin{align}
J_3
&=J_3(p_0,p_1,p_2,p_3;q_1,q_2,q_3)
\nonumber\\
&=\iiint\limits_{[0,1]^3}
\frac{y_1^{p_1}(1-y_1)^{q_1}y_2^{p_2}(1-y_2)^{q_2}y_3^{p_3}(1-y_3)^{q_3}\,\d y_1\,\d y_2\,\d y_3}
{(1-y_3(1-y_1y_2))^{p_0+1}}
\label{J3}
\end{align}
is the generalised Beukers integral \cite{Be79}. As shown in \cite{RV01}, the condition
\begin{equation}
p_3+q_3=q_1+q_2
\label{cond1}
\end{equation}
on its parameters implies that $J_3=2A\zeta(3)-B$; furthermore, by the results in~\cite{Zu04} we have
\begin{align*}
A
&=A(p_0,p_1,p_2,p_3;q_1,q_2,q_3)
\\
&=(-1)^{p_0+p_1+p_2+p_3}\sum_{k\in\mathbb Z}\binom{k}{p_0}\binom{k+q_3-p_0}{p_3+q_3-p_0}\binom{q_1}{k-p_1}\binom{q_2}{k-p_2},
%\label{sumA}
\end{align*}
which combined with \eqref{I3} leads to another proof of formula \eqref{sumQ}.

There are several other important consequences of the explicit reduction of $I''(\ba)$ to $J_3$ given in~\eqref{I3}.
The formula provides us with access to the arithmetic of $I''(\ba)$ (for example, it `explains' the appearance of the $d_{2n}$ factor in~\eqref{dn-totsym}). 
It also means that any $\mathbb Q$-linear relations between the integrals $I(\ba)$ deduced from manipulating  the integrand (for example, contiguity relations) also hold for $I'(\ba)$, $I''(\ba)$ and the rational coefficients $Q(\ba)$, $\hat P(\ba)$, $P(\ba)$. 

It is important to point out that the integrals \eqref{J3}, subject to \eqref{cond1},  are equivalent to  generalised cellular integrals on $\mathcal{M}_{0,6}$ \cite{Br16}.   The formulae described in this section are indicative of a general recursive structure between cellular integrals on $\mathcal{M}_{0,n}$ for different $n$, which should reflect  the recursive structure of their boundary strata.

\section({Group structure for \003\266(5)})%
{Group structure for $\zeta(5)$}
\label{group}

In addition to the automorphism $\fri_1$ recorded in Section~\ref{sec3},  the families $I(\ba)$ and $J(\bp;\bq)$    admit `hypergeometric' transformations. 
These relate integrals within a given family up  to a rational  prefactor which may be expressed as  a quotient of factorials.
Such extended automorphisms play an important role in extracting  arithmetic information about the coefficients $Q,\hat P,P$ in the decomposition \eqref{deco} of $I(\ba)$.

First observe that the $_3F_2$-hypergeometric functions in representation \eqref{3F2J} of $J(\bp;\bq)$
are symmetric with respect to  permutations of their top parameters.
Such permutations however affect the factor $p_1!\,q_1!\,p_2!\,p_4!\,q_4!\,p_5!\,q_5!$\,.

A simple way to describe these permutations is through manipulations of the internal double integral over $y_1,y_2$ in~\eqref{J1}:
\begin{align*}
&
\iint\limits_{[0,1]^2}\frac{y_1^{p_1}(1-y_1)^{q_1}y_2^{p_2}(1-y_2)^{q_2}}
{(1-y_3(1-y_1y_2))^{p_0+1}}\,\d y_1\,\d y_2
\\ &\quad
=\frac1{(1-y_3)^{p_0+1}}
\,\frac{p_1!\,p_2!\,q_1!\,q_2!}{(p_1+q_1+1)!\,(p_2+q_2+1)!}
\,{}_3F_2\bigg(\begin{matrix} p_0+1, \, p_1+1, \, p_2+1 \\ p_1+q_1+2, \, p_2+q_2+2 \end{matrix}\biggm|\frac{-y_3}{1-y_3}\bigg).
\end{align*}
The symmetry of the parameters $p_1+1$ and $p_2+1$ in the $_3F_2$-representation implies that 
\begin{align*}
&
\iint\limits_{[0,1]^2}\frac{y_1^{p_1}(1-y_1)^{q_1}y_2^{p_2}(1-y_2)^{q_2}}
{(1-y_3(1-y_1y_2))^{p_0+1}}\,\d y_1\,\d y_2
\\ &\quad
=\frac{q_1!\,q_2!}{(p_1-p_2+q_1)!\,(p_2-p_1+q_2)!}
\iint\limits_{[0,1]^2}\frac{x_1^{p_2}(1-x_1)^{p_1-p_2+q_1}x_2^{p_1}(1-x_2)^{p_2-p_1+q_2}}
{(1-y_3(1-y_1y_2))^{p_0+1}}\,\d y_1\,\d y_2;
\end{align*}
the symmetry of the parameters $p_0+1$ and $p_1+1$ leads to
\begin{align*}
&
\iint\limits_{[0,1]^2}\frac{y_1^{p_1}(1-y_1)^{q_1}y_2^{p_2}(1-y_2)^{q_2}\,\d y_1\,\d y_2}
{(1-y_3(1-y_1y_2))^{p_0+1}}
\\ &\quad
=(1-y_3)^{p_1-p_0}
\,\frac{p_1!\,q_1!}{p_0!\,(p_1-p_0+q_1)!}
\iint\limits_{[0,1]^2}\frac{x_1^{p_0}(1-x_1)^{p_1-p_0+q_1}x_2^{p_2}(1-x_2)^{q_2}}
{(1-y_3(1-y_1y_2))^{p_1+1}}\,\d y_1\,\d y_2.
\end{align*}
Inserting these findings in the 5-fold integral \eqref{J1} we see that the quantity
\[
\frac{J(\bp;\bq)}{p_1!\,p_2!\,q_1!\,q_2!}
\]
is invariant under the group (of order~6) generated by two involutions
\[
\frp_{12}\colon(\bp;\bq)
\mapsto(p_0,p_2,p_1,p_3,p_4,p_5,p_6;p_1+q_1-p_2,p_2+q_2-p_1,q_3,q_4,q_5)
\]
and
\[
\frp_{01}\colon(\bp;\bq)
\mapsto(p_1,p_0,p_2,p_3,p_4,p_5,p_6;p_1+q_1-p_0,q_2,p_1+q_3-p_0,q_4,q_5).
\]
These transformations respect both conditions \eqref{cond12} and \eqref{cons1}.
With the reflection \eqref{i1pq} we also get another group of order~6 generated by two involutions
\[
\frp_{45}=\fri_1\frp_{12}\fri_1\colon(\bp;\bq)
\mapsto(p_0,p_1,p_2,p_3,p_5,p_4,p_6;q_1,q_2,q_3,p_4+q_4-p_5,p_5+q_5-p_4)
\]
and
\[
\frp_{56}=\fri_1\frp_{01}\fri_1\colon(\bp;\bq)
\mapsto(p_0,p_1,p_2,p_3,p_4,p_6,p_5;q_1,q_2,p_5+q_3-p_6,q_4,p_5+q_5-p_6).
\]
All these transformations are in the group $\langle\fri_1,\frp_{01},\frp_{12}\rangle$ of order $2\times3!\times3!=72$; they respect both \eqref{cond12} and \eqref{cons1} and keep the quantity
\[
\frac{J(\bp;\bq)}{p_1!\,p_2!\,p_4!\,p_5!\,q_1!\,q_2!\,q_4\,q_5!}
\]
invariant. This group can be regarded as a `trivial' hypergeometric group as it only takes into account trivial symmetries of the underlying hypergeometric representation. It acts on both the 12-parameter family $J(\bp;\bq)$ and 8-parameter family $I(\ba)$.

There are more transformations available for the 8-parameter family $I(\ba)$, although they can be observed more easily on the corresponding subfamily of $J(\bp;\bq)$.
Conditions~\eqref{cond12} imply that both integrals
\begin{equation*}
\frac1{2\pi i}\,\raisebox{1.8mm}{$\displaystyle\bigintsss_{-c_1-i\infty}^{-c_1+i\infty}$}
\frac{\splitfrac{\Gamma(p_0+1+s)\Gamma(p_1+1+s)\Gamma(p_2+1+s)\Gamma(p_3+2+t+s)}{\times\Gamma(q_3-p_0-p_6-1-t-s)\Gamma(-s)}}
{\Gamma(p_1+q_1+2+s)\Gamma(p_2+q_2+2+s)}\,\d s
%\label{Is}
\end{equation*}
and
\[
\frac1{2\pi i}\,\raisebox{1.8mm}{$\displaystyle\bigintsss_{-c_2-i\infty}^{-c_2+i\infty}$}
\frac{\splitfrac{\Gamma(p_4+1+t)\Gamma(p_5+1+t)\Gamma(p_6+1+t)\Gamma(p_3+2+s+t)}{\times\Gamma(q_3-p_0-p_6-1-s-t)\Gamma(-t)}}
{\Gamma(p_4+q_4+2+t)\Gamma(p_5+q_5+2+t)}\,\d t
\]
in the double integration in~\eqref{intJ} are subject to Bailey's transformation \cite[\S\,6.3, eq.~(2)]{Ba35} (see also \cite[eq.~(4.7.1.3)]{Sl66}) and therefore can be cast as $_7F_6$ very well-poised hypergeometric series:
\begin{align}
&
\frac1{2\pi i}\,\raisebox{1.8mm}{$\displaystyle\bigintsss_{-i\infty}^{i\infty}$}
\frac{\splitfrac{\Gamma(b+t)\,\Gamma(c+t)\,\Gamma(d+t)\,\Gamma(1+a-e-f+t)}{\times\Gamma(1+a-b-c-d-t)\,\Gamma(-t)\,\d t}}
{\Gamma(1+a-e+t)\,\Gamma(1+a-f+t)}
\nonumber\\ &\;
=\frac{\splitfrac{\Gamma(1+a)\,\Gamma(b)\,\Gamma(c)\,\Gamma(d)\,\Gamma(1+a-c-d)\,\Gamma(1+a-b-d)}{\times\Gamma(1+a-b-c)\,\Gamma(1+a-e-f)}}
{\Gamma(1+a-b)\,\Gamma(1+a-c)\,\Gamma(1+a-d)\,\Gamma(1+a-e)\,\Gamma(1+a-f)}
\nonumber\\ &\;\quad\times
{}_7F_6\biggl(\begin{alignedat}{7} & a, &\,& 1+\tfrac12a, &\,& \quad\; b, &\,& \qquad c, &\,& \qquad d, &\,& \qquad e, &\,& \qquad f \\[-3\jot]
 & &\,& \quad \tfrac12a, &\,& {\kern-2.5mm}1+a-b, &\,& 1+a-c, &\,& 1+a-d, &\,& 1+a-e, &\,& 1+a-f \end{alignedat} \biggm| 1 \biggr).
\label{7F6tr}
\end{align}
Note that it is this transformation which plays a crucial role in a hypergeometric interpretation \cite{Zu04} of the Rhin--Viola group \cite{RV01} for the linear forms in $1$ and $\zeta(3)$.

If we choose
\begin{gather*}
a=p_2+q_1+q_2+2, \;\; b=p_0+1, \;\; c=p_3+t+2, \;\; d=p_2+1, \\ e=p_2-p_1+q_2+1, \;\; f=q_1+1
\end{gather*}
in \eqref{7F6tr}, and then apply the same identity with $d$ and $e$ interchanged\,---\,this does not affect the $_7F_6$ series\,---\,we find out that
\begin{align*}
&
\frac1{2\pi i}\,\raisebox{1.8mm}{$\displaystyle\bigintsss_{-c_1-i\infty}^{-c_1+i\infty}$}
\frac{\splitfrac{\Gamma(p_0+1+s)\Gamma(p_1+1+s)\Gamma(p_2+1+s)}{\times\Gamma(p_3+2+t+s)\Gamma(q_3-p_0-p_6-1-t-s)\Gamma(-s)}}
{\Gamma(p_1+q_1+2+s)\Gamma(p_2+q_2+2+s)}\,\d s
\\ &\;
=\frac{p_1!\,p_2!\,(q_1+q_2-p_0)!\,\Gamma(q_1+q_2-p_3-t)}
{q_2!\,(p_1+q_1-p_0)!\,(p_2+q_2-p_1)!\,\Gamma(p_1+q_1-p_3-t)}
\\ &\;\quad\times
\frac1{2\pi i}\,\raisebox{1.8mm}{$\displaystyle\bigintsss_{-c_1'-i\infty}^{-c_1'+i\infty}$}
\frac{\splitfrac{\Gamma(p_0+1+s)\Gamma(q_2+1+s)\Gamma(p_2+q_2-p_1+1+s)}{\times\Gamma(p_3+2+t+s)\Gamma(p_1-p_0+q_1-p_3-1-t-s)\Gamma(-s)}}
{\Gamma(q_1+q_2+2+s)\Gamma(p_2+q_2+2+s)}\,\d s.
\end{align*}
Assuming now \eqref{cons1}, that is, using $p_1+q_1=p_3$, and substituting into the double integral expression for $J(\bp;\bq)$ we obtain
\begin{align*}
J(\bp;\bq)
&=\frac{p_1!\,p_2!\,q_1!\,q_4!\,q_5!}{p_0!\,p_6!\,(p_3-p_0)!\,(p_2+q_2-p_1)!}
\\ &\;\times
\frac1{2\pi i}\int_{-c_1'-i\infty}^{-c_1'+i\infty}\d s
\,\frac{\Gamma(p_0+1+s)\Gamma(q_2+1+s)\Gamma(p_2+q_2-p_1+1+s)\Gamma(-s)}
{\Gamma(q_1+q_2+2+s)\Gamma(p_2+q_2+2+s)}
\\ &\;\quad\times
\frac1{2\pi i}\int_{-c_2-i\infty}^{-c_2+i\infty}\d t
\,\frac{\Gamma(p_4+1+t)\Gamma(p_5+1+t)\Gamma(p_6+1+t)\Gamma(q_1+q_2-p_3-t)}{\Gamma(p_4+q_4+2+t)\Gamma(p_5+q_5+2+t)}
\\ &\;\qquad\times
\Gamma(p_3+2+s+t)\Gamma(-p_0-1-t-s)
\\ \intertext{(after the shift $t\mapsto t+r$ with $r=q_1+q_2-p_3=q_3-p_6$)}
&=\frac{p_1!\,p_2!\,q_1!\,q_4!\,q_5!}{p_0!\,p_6!\,(p_3-p_0)!\,(p_2+q_2-p_1)!}
\\ &\;\times
\frac1{2\pi i}\int_{-c_1'-i\infty}^{-c_1'+i\infty}\d s
\,\frac{\Gamma(p_0+1+s)\Gamma(q_2+1+s)\Gamma(p_2+q_2-p_1+1+s)\Gamma(-s)}
{\Gamma(q_1+q_2+2+s)\Gamma(p_2+q_2+2+s)}
\\ &\;\quad\times
\frac1{2\pi i}\int_{-c_2'-i\infty}^{-c_2'+i\infty}\d t
\,\frac{\Gamma(p_4+r+1+t)\Gamma(p_5+r+1+t)\Gamma(p_6+r+1+t)\Gamma(-t)}{\Gamma(p_4+q_4+r+2+t)\Gamma(p_5+q_5+r+2+t)}
\\ &\;\qquad\times
\Gamma(p_3+r+2+s+t)\Gamma(-p_0-r-1-t-s)
\displaybreak[2]\\
&=\frac{p_1!\,p_2!\,q_1!\,q_4!\,q_5!}{p_0!\,p_6!\,(p_3-p_0)!\,(p_2+q_2-p_1)!}
\cdot\frac{p_0'!\,p_6'!\,(q_1'+q_2'-p_0')!}{q_1'!\,q_2'!\,q_4'!\,q_5'!}
\cdot J(\bp';\bq')
\\
&=\frac{p_2!\,q_3!}{p_6!\,(p_2+q_2-p_1)!}
\cdot J(\bp';\bq'),
\end{align*}
where
\begin{align*}
(\bp';\bq')
&=(p_0,q_2,p_2+q_2-p_1,p_3+r,p_4+r,p_5+r,p_6+r; q_1,p_1,q_3-r,q_4,q_5)
\\
&=(p_0,q_2,p_2+q_2-p_1,q_1+q_2,p_4-p_6+q_3,p_5-p_6+q_3,q_3; q_1,p_1,p_6,q_4,q_5).
\end{align*}
The identity can be interpreted as the invariance of
\[
\frac{J(\bp;\bq)}{p_2!\,q_3!}
\]
under the involution $\frh'\colon(\bp;\bq)\mapsto(\bp';\bq')$.
A different manipulation with Bailey's transformation \eqref{7F6tr}, more in line with the normalisation in \cite[Sect.~4]{Zu04}, leads to a different hypergeometric involution
\begin{align*}
\frh\colon(\bp;\bq)
&\mapsto(q_2,p_1,p_2-p_0+q_2,p_3-p_0+q_2,p_4-p_0+q_2,p_5-p_0+q_2,p_6-p_0+q_2;
\\ &\qquad
q_1-p_0+q_2,p_0,q_3,q_4,q_5)
\end{align*}
which keeps the quantity
\[
\frac{J(\bp;\bq)}{p_2!\,q_1!\,q_2!\,(p_6-p_0+q_2)!}
\]
invariant.

The group $\frG$ generated by $\fri_1$, $\frp_{01}$, $\frp_{12}$ and $\frh$ (or by $\frh'$ instead of the latter) acts naturally on the 8-parameter set $(\bp;\bq)$ subject to the constraints \eqref{cond12}, \eqref{cons1}.
Because the 8-parameter set is in one-to-one correspondence with the original parameter set $\ba=(a_1,\dots,a_8)$, the action of $\frG$ can be translated into one on the parameters~$\ba$.
We obtain
\begin{align*}
\fri_1\colon\ba&\mapsto(a_5,a_4,a_3,a_2,a_1,a_7,a_6,-a_2+a_4-a_6+a_7+a_8),
\\
\frp_{01}\colon\ba&\mapsto(a_1,a_2,-a_2+a_4+a_5,a_2+a_3-a_5,a_5,a_6,a_7,-a_2-a_3+a_4+a_5+a_8),
\\
\frp_{12}\colon\ba&\mapsto(a_1,a_2,-a_2+a_4+a_8,a_2+a_3-a_8,-a_2-a_3+a_4+a_5+a_8,
\\ &\qquad
a_2+a_3-a_4+a_6-a_8,a_7,a_8),
\\
\frh\colon\ba&\mapsto(a_1,a_2,a_3+a_6-a_8,a_4-a_6+a_8,a_5+a_6-a_8,a_8,-a_6+a_7+a_8,a_6),
\\
\frh'\colon\ba&\mapsto(a_1,a_2,a_3,a_4,a_2+a_3-a_4+a_6-a_8,
\\ &\qquad
-a_2-a_3+a_4+a_5+a_8,-a_2-a_3+a_4+a_5-a_6+a_7+a_8,a_8).
\end{align*}
It can be checked that the group is a permutation group on the 28~elements
\begin{equation}
\begin{gathered}
h_i=a_i \quad\text{for}\; i=1,\dots,8,
\\
\begin{alignedat}{3}
h_9&=a_1+a_2-a_4, &\quad
h_{10}&=a_1+a_5-a_3, &\quad
h_{11}&=a_1+a_8-a_3, \\
h_{12}&=a_2+a_3-a_5, &\quad
h_{13}&=a_2+a_3-a_8, &\quad
h_{14}&=a_3+a_6-a_8, \\
h_{15}&=a_3+a_4-a_1, &\quad
h_{16}&=a_4+a_5-a_2, &\quad
h_{17}&=a_4+a_8-a_6, \\
h_{18}&=a_4+a_8-a_2, &\quad
h_{19}&=a_5+a_6-a_8, &\quad
h_{20}&=a_7+a_8-a_6,
\end{alignedat}
\\
\begin{alignedat}{2}
h_{21}&=a_1+a_2+a_6-a_4-a_8, &\quad
h_{22}&=a_1+a_7+a_8-a_3-a_6, \\
h_{23}&=a_2+a_3+a_6-a_4-a_8, &\quad
h_{24}&=a_2+a_3+a_6-a_7-a_8, \\
h_{25}&=a_4+a_5+a_8-a_2-a_3, &\quad
h_{26}&=a_4+a_7+a_8-a_2-a_6, \\
h_{27}&=a_4+a_7+2a_8-a_2-a_3-a_6, &\quad
h_{28}&=a_4+a_5+a_7+a_8-a_2-a_3-a_6.
\end{alignedat}
\end{gathered}
\label{eq:h}
\end{equation}
The hyperplanes $h_j=0$ for $j=1,\dots,28$ are precisely orbits of $a_i=0$ for $i=1,\dots,8$ under the action of~$\frG$. There is one additional  hyperplane
\[
a_1+a_5+a_7+a_8-a_2-a_3=0
\]
which is also 
 preserved by this action.
The analysis above implies that the quantity
\begin{equation}
\frac{I(\ba)}{\prod_{i\in\cF}h_i!},
\quad\text{where}\;\;
\cF=\{1,2,3,4,5,6,7,9,10,11,14,16,18,20,23,27,28\},
\label{F17}
\end{equation}
is invariant under $\frG$.

\begin{remark}
\label{rem17}
Notice that the convergence conditions on the integral $I(\ba)$ in the introduction, namely, the non-negativity of the seventeen quantities~\eqref{17elem}, are precisely $h_i\ge0$ for $i\in\cF$.
\end{remark}

This invariance can alternatively be stated for the quantity
\[
%\frac{J(\bp;\bq)}{\splitfrac{p_1!\,p_2!\,p_4!\,p_5!\,q_1!\,q_2!\,q_3!\,q_4!\,q_5!\,(p_4-p_0+q_2)!\,(p_2-p_6+q_4)!\,(p_2-p_0+q_3)!\,(p_4-p_6+q_3)!}
%{\times(p_2+p_4-p_0-p_6+q_3)!\,(q_2-q_3+q_4)!\,(q_2-q_5+q_1)!\,(q_4-q_1+q_5)!}}
\frac{J(\bp;\bq)}{\begin{aligned}
& p_1!\,p_2!\,p_4!\,p_5!\,q_1!\,q_2!\,q_3!\,q_4!\,q_5!\,(p_4-p_0+q_2)!\,(p_2-p_6+q_4)!\,
\\[-1.5mm] &\qquad\times (p_2-p_0+q_3)!\,(p_4-p_6+q_3)!\,(p_2+p_4-p_0-p_6+q_3)!\,
\\[-1.5mm] &\qquad\times (q_2-q_3+q_4)!\,(q_2-q_5+q_1)!\,(q_4-q_1+q_5)!
\end{aligned}}
\]
provided that the parameters $(\bp;\bq)$ are subject to the conditions~\eqref{cond12}, \eqref{cons1}.
Furthermore, the 28-element multiset $\bh$ naturally splits into two submultisets
\begin{align*}
\bh'=\frG a_1
&=\{h_1,h_3,h_5,h_6,h_7,h_8,h_9,h_{10},h_{11},h_{14},h_{16},h_{18},
\\ &\qquad\qquad
h_{19},h_{20},h_{21},h_{22},h_{23},h_{25},h_{26},h_{27},h_{28}\}
\end{align*}
of size 21 and
\[
\bh''=\frG a_2
=\{h_2,h_4,h_{12},h_{13},h_{15},h_{17},h_{24}\}
\]
of size 7,  where $\bh'\cap\bh''=\emptyset$. Each of these multisets  is individually acted  upon by the above symmetry group $\frG$. Also notice that
\[
\{h_i:i\in\cF\}\cap\bh'=\{h_1,h_3,h_5,h_6,h_7,h_9,h_{10},h_{11},h_{14},h_{16},h_{18},h_{20},h_{23},h_{27},h_{28}\}
\]
and
\[
\{h_i:i\in\cF\}\cap\bh''=\{h_2,h_4\}
\]
have size 15 and 2, respectively (while their complements in $\bh'$ and $\bh''$ have size 6 and~5).

\begin{remark} 

\label{rem-num}
It is a good moment to bring some  related numerology to the reader's attention. 
Firstly, we  recognise the sizes 28, 21 and 7 of the multisets $\bh$, $\bh'$ and $\bh''$ as being $\binom82$, $\binom72$ and $\binom71$, respectively.
Secondly, the order of the group $\frG$ is $7!=5040$, which is  the order of the symmetric group $\Sigma_7$ on 7 letters. The precise explanation for this will be made apparent in the next sections, as we shall show that indeed $\frG$ is naturally isomorphic to $\Sigma_7$. 
Note that the group $\Sigma_7$ is also isomorphic   to  the Weyl group $W(\mathrm A_6)$ of the root system $\mathrm A_6$, which is consistent with the known hypergeometric groups for rational approximations to $\zeta(2)$ (as in \cite{RV96}), to $\zeta(3)$ (as in \cite{RV01}) and to $\zeta(4)$ (as in \cite{MZ20}): 
they are equal to $|W(\mathrm A_4)|=120$, $|W(\mathrm D_5)|=1920$ and $|W(\mathrm E_6)|=51840$, respectively.
Precisely how this pattern of symmetry groups extends to more general cellular integrals (as we expect) remains a mystery. 
\end{remark}

Now recall that by the results in Section~\ref{zeta(3)} the group $\frG$ acts not only on the integrals $I(\ba)$ but also on their parts $I'(\ba)\in\mathbb Q\zeta(5)+\mathbb Q$ and $I''(\ba)\in\mathbb Q\zeta(3)+\mathbb Q$, where
$I(\ba)=2I'(\ba)+4I''(\ba)\zeta(2)$.
As discussed in Remark~\ref{rem-decom},  it is in principle possible to prove, with considerable effort, the experimental observation:
%\francisnote{it is not clear if we are claiming that this is rigourously proven, or merely achievable with some considerable effort?} \wadimnote{It is achievable with some serious effort... when it is worth that and I am retired to spend years on that.}
\begin{equation}
d_{m_1n}d_{m_2n}d_{m_3n}d_{m_4n}d_{m_5n}I'(\ba\cdot n)\in\mathbb Z\zeta(5)+\mathbb Z,
\label{incl}
\end{equation}
where $m_1\ge m_2\ge m_3\ge m_4\ge m_5$ are five consecutive maxima of the corresponding 28-element multiset $h_1,\dots,h_{28}$ and $d_N$ denotes the least common multiple of $1,2,\dots,N$.
Because the latter multiset is invariant under $\frG$, we conclude from \eqref{incl} that
\begin{equation*}
d_{m_1n}d_{m_2n}d_{m_3n}d_{m_4n}d_{m_5n}I'(\frg\ba\cdot n)\in\mathbb Z\zeta(5)+\mathbb Z
\quad\text{for any}\; \frg\in\frG;
\end{equation*}
in particular,
\begin{align*}
&
\frac{\prod_{i\in\cF}(\frg h_i)!}{\prod_{i\in\cF}h_i!}
\cdot d_{m_1n}d_{m_2n}d_{m_3n}d_{m_4n}d_{m_5n}I'(\ba\cdot n)
\\ &\qquad\qquad
=d_{m_1n}d_{m_2n}d_{m_3n}d_{m_4n}d_{m_5n}I'(\frg\ba\cdot n)\in\mathbb Z\zeta(5)+\mathbb Z
\quad\text{for any}\; \frg\in\frG.
\end{align*}
This means that if we take
\begin{equation}
\nu_p=\max_{\frg\in\frG}\ord_p\frac{\prod_{i\in\cF}h_i!}{\prod_{i\in\cF}(\frg h_i)!}
\label{nu_p}
\end{equation}
and consider the quantity
\[
\Phi_n=\Phi_n(\ba)=\prod_{p>\sqrt{m_1n}}p^{\nu_p},
\]
then
\begin{equation}
\Phi_n^{-1}\cdot d_{m_1n}d_{m_2n}d_{m_3n}d_{m_4n}d_{m_5n}I'(\ba\cdot n)
\in\mathbb Z\zeta(5)+\mathbb Z.
\label{sharp_incl}
\end{equation}
Furthermore, we can compute the limit of $|I'(\ba\cdot n)|^{1/n}$ as $n\to\infty$ using the asymptotics of $J(\bp;\bq)$ discussed in Section~\ref{sec4}.

\begin{remark}
\label{m_i}
The above choice of $m_1,\dots,m_5$ may be not optimal. It is plausible to expect that, for an appropriate value of  $\ell\in\{1,\dots,5\}$, one can take it as follows: $m_1\ge\dots\ge m_\ell$ are successive maxima of the 21-element multiset $\bh'$ and $m_{\ell+1}\ge\dots\ge m_5$ are successive maxima of the 7-element multiset $\bh''$.
Notice that  for any such $\ell$, the set  $\{m_1,\dots,m_5\}$ is  preserved under the action of~$\frG$. A similar situation arises in the context of  multi-parameter rational approximations to $\zeta(2)$ \cite{RV96,Zu14}, $\zeta(3)$ \cite{RV01} and $\zeta(4)$ \cite{MZ20}.
\end{remark}

\section{Automorphisms of the asymptotic polynomials}
\label{sec:Graph}
Consider the  discriminant of  the cubic polynomial defined in \eqref{eqn: Fcubic}, viewed as a  polynomial function of  the parameters $a_1,\dots, a_8$. Upon dividing out (the squares of) linear factors of the form $a^2_5, a^2_8,(a_2+a_3-a_4-a_5-a_8)^2$, one is left with an  irreducible polynomial $D\in \mathbb{Q}[a_0,\dots, a_8]$  which is homogeneous of degree 16. We wish to  compute the subgroup $\operatorname{Aut}_D \leq \operatorname{GL}_8(\mathbb{Z})$ of linear automorphisms  acting on the parameters $a_0,\dots,a_8$ which preserves  the equation $D=0$. 

To this end, we  consider hyperplanes $H$, defined by a  homogeneous linear form in the parameters $a_1,\dots, a_8$ with the property that the restriction of the discriminant $D\big|_H$ 
to $H$ factorizes  into a product  $D_6  D_5^2$, where  $D_6, D_5$ are irreducible of degree $6$ and~$5$. 
We find that the hyperplanes of this form are  the 28 hyperplanes $h_i=0$ with $h_1,\dots, h_{28}$ listed in \eqref{eq:h}, as well as the exceptional hyperplane 
\[ a_5 +a_7 +a_8 +a_1-a_2 -a_3 =0. \]
Since the action of $\operatorname{Aut}_D$ is linear, it  permutes this set of  29 hyperplanes. 
In order to understand this action, we consider the restriction $D|_{H_i\cap H_j}$ to all pairs of such hyperplanes $H_i,H_j$. We find that for certain pairs $H_i, H_j$ (which we call of the `first kind'), we have  a factorisation of the form
\[ D\big|_{H_i\cap H_j} = (P_1P'_1P_3P'_3)^2\]
where $P_1,P'_1$ are linear and $P_3,P'_3$ are of degree~$3$. For other pairs, another type of factorisation may also occur, whereby $D|_{H_i\cap H_j}$ has 5 linear factors of multiplicity two, and a single irreducible component of degree $6$. A third type of factorisation  occurs if and only if  one of  $H_i, H_j$ is the exceptional hyperplane.

\begin{figure}[ht]
%\quad {\includegraphics[width=12cm]{Graph168.png}} 
%\quad {\includegraphics[width=8cm]{BlueGraph168.png}} 
\includegraphics[width=8cm]{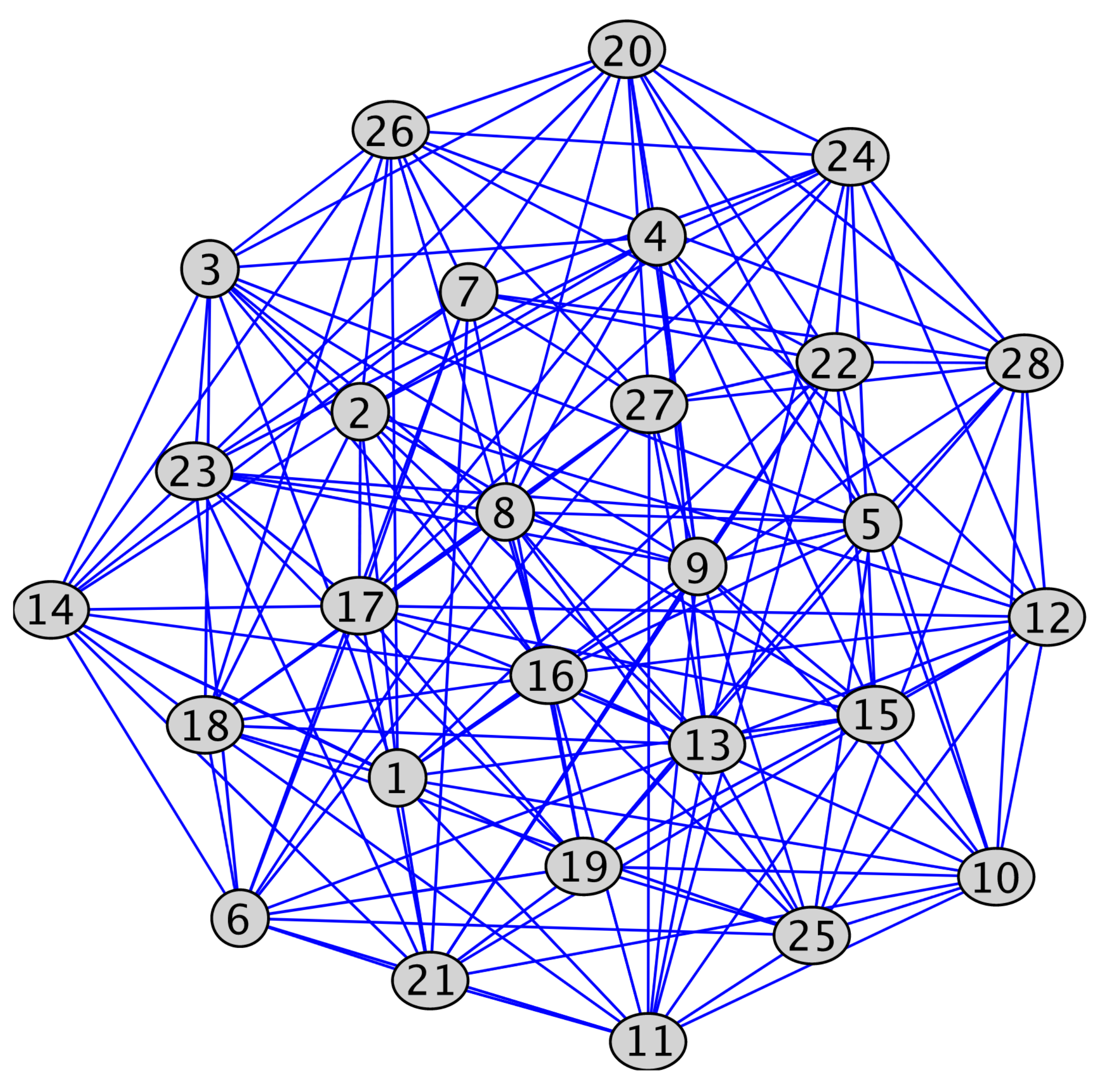}
\caption{The graph $G$.}
\label{Graph168}
\end{figure}

The action of $\operatorname{Aut}_D$ thus preserves the zero loci of the set of linear forms  $h_1,\dots, h_{28}$ together with the data  of which pairs $(h_i,h_j)$  are of the first kind. 
This data may be encoded by a graph $G$, with one vertex for every hyperplane $H_i = V(h_i)$, and with an edge between vertices $H_i$ 
and $H_j$ if and only if $(h_i,h_j)$ are of the first kind. It is depicted in Figure \ref{Graph168}, where the vertex  $H_i$ is denoted by $i$.

\begin{figure}[ht]
\quad {\includegraphics[width=12cm]{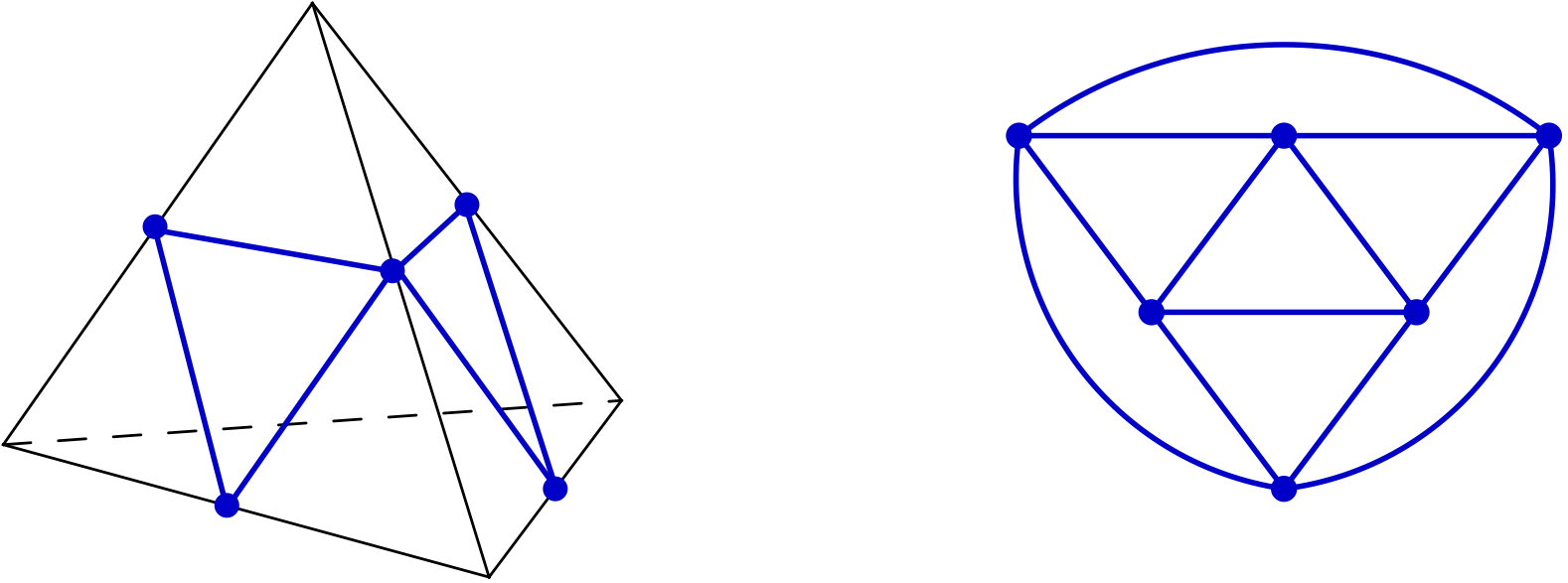}} 
\put(-350,20){$s_1$}\put(-245,-10){$s_2$}\put(-200,40){$s_3$}\put(-280,135){$s_0$}
\caption{The complete graph $K_4$ on 4 vertices labelled $s_0,\ldots, s_3$ is depicted on the left as a tetrahedron. The two front faces have been subdivided into triangles. Subdividing all faces leads to the graph $G_4$  with 6 vertices and 12 edges shown on  the right. }
\label{figG_n}
\end{figure}
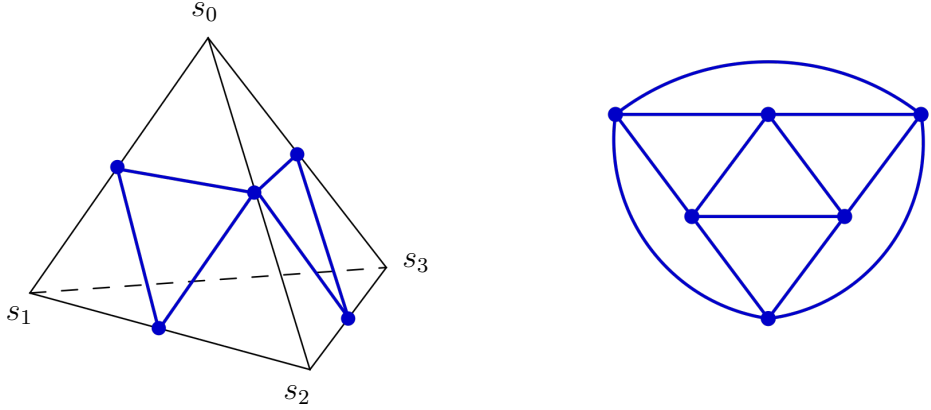

One finds that this graph  has 168 edges and that  each vertex has degree exactly~12. Furthermore, one can easily check by computer   that  its automorphism group is of order $8!$\,.  On observing that $168 = 3 \binom{8}{3}$ and $28 =\binom{8}{2}$, one is led to suspect that $G$ may be isomorphic to the highly symmetric graph $G_8$ defined as follows:

Let  $K_n$ denote the complete graph with $n\geq 2$ vertices. Define a new graph $G_n$
 whose vertices are given by the set of edges  $\{i,j\}$ of $K_n$, with an edge between $\{i,j\}$ and $\{k,\ell\}$  if and only if $\{i,j,k,\ell\}$ is a set with 3 elements. Equivalently, the graph $G_n$ may be constructed by placing a single vertex  along every edge of the complete graph $K_n$, and connecting every pair of such vertices which lie on the same face of~$K_n$. This construction is illustrated in Figure~\ref{figG_n}.

Once one suspects that $G$ might be isomorphic to $G_8$, it takes only a  little  detective work to find an isomorphism.
Indeed, if we label the vertices of $G_8$ by $s_0,s_1,\dots,s_7$ we find that an explicit isomorphism $G \cong G_8$ is given by the following table:
\[ 
\begin{array}{c|cccccccc}
  & s_0 & s_1 & s_2 & s_3 & s_4 & s_5 & s_6 & s_7  \\ \hline
  s_0    & &  15  &  2  & 4  & 12  & 13  & 17    & 24          \\
 s_1     & &  &  1     & 9  & 10  & 11  & 21    & 22    \\
 s_2     & &     &  & 3  & 16  & 18  & 14    & 26    \\ 
  s_3    & &    &   &   & 5  & 8  & 23    & 20    \\
  s_4    & &    &   &   &   & 25 & 19    & 28    \\
  s_5    & &    &   &   &   &  & 6    & 27    \\
  s_6    & &    &   &   &   &  &     & 7    
 \end{array}
\]
The table means, for example,  that the  vertex of $G$ assigned to hyperplane $h_1$ corresponds to the vertex $\{s_1,s_2\}$ of $G_8$, and the vertex assigned to $h_{27}$ corresponds to $\{s_5,s_7\}$.

From this it is clear that  the automorphism group $\operatorname{Aut}(G) \cong \operatorname{Aut}(G_8)$ is  the symmetric group of order $8!$ acting via permutations on $S=\{s_0,\ldots, s_7\}$. Furthermore, we discover  that the hypergeometric group $\frG$ is  precisely the subgroup of order $7!$ which fixes the vertex $s_0$.
In fact, we find  with a little more calculation that 
\[  \frG =  \operatorname{Aut}(G) \cap \operatorname{GL}_8(\mathbb{Z})^+\]
where $\operatorname{GL}_8(\mathbb{Z})^+$ denotes the subgroup of linear transformations $g \in \operatorname{GL}_8(\mathbb{Z})$  which are positive.%
\footnote{In the symmetric coordinates $s_i$,  an element is positive if an only if it fixes $s_0$; see \eqref{H-matr}.}
This means the following: a linear form $\ell(a_1,\dots, a_8)$ is called positive if in the symmetric case $a_1=a_2=\dots= a_8=n$, where $n>0$, it satisfies $\ell(n,\ldots, n)>\nobreak0$. An element $g$ is positive if it sends a positive defining equation of every $h_i$ to a positive linear form. The above calculation implies that  the group $\frG$ is largest possible and that we have indeed exhausted all  symmetries of \eqref{I1}. 
   
  Nevertheless, there exist elements in $\operatorname{Aut}_D$ which are not positive. For example, by considering a permutation which does not fix $s_0$,  we are led to the  automorphism 
   \[ (a_1,a_2,a_4) \mapsto (-a_2, -a_1, a_4-a_1-a_2  ), \]
   which fixes all other $a_i$,  preserves the exceptional hyperplane,   acts upon the set of hyperplanes $h_i=0$, for $1\leq i \leq 28$, and preserves the discriminant locus $D=0$.  It is very unclear whether the existence of this curious automorphism  has any implications for the arithmetic of linear forms in $\zeta(5)$. In any case, we check that this non-positive automorphism, together with the group $\frG$, generates the full group of symmetries of the graph $G$ of order~$8!$\,. 

\section{Duality}
\label{sec-duality}

A general duality  was proven for \emph{totally symmetric} cellular integrals in \cite{Br16}. A natural definition of duality for generalised cellular integrals is as follows.
 
 Let $N\geq 5$ and let $\pi$ be a convergent permutation on $\{1,\dots,N\}$ in the sense of \cite{Br16}.%
\footnote{A permutation $\pi$ is convergent if and only if the corresponding cellular integral is finite. This can equivalently be formulated as a combinatorial condition on the permutation: roughly speaking,  consecutive sequences of integers $n,n+1,\ldots, n+k$ modulo $N$ must be  non-consecutive after applying $\pi$, for all  $1\leq k\leq N-2$.} 
In that paper,  the generalised  rational  cellular function was defined by
\begin{equation*}
%\label{Genf}
f_{\pi}(\ba,\bb) = \prod_{i\in \Z/N\Z} \frac{ (z_i - z_{i+1})^{a_{i,i+1}}}{ (z_{\pi_i} - z_{\pi_{i+1}})^{b_{\pi_i,\pi_{i+1}}}}
\end{equation*}
where the multi-indices $\ba$, $\bb$ are subject to homogeneity equations 
\[ a_{i-1,i}  + a_{i,i+1} = b_{\pi_{j-1}, \pi_{j}} + b_{\pi_{j}, \pi_{j+1}}\]
for all $i,j$ such that  $\pi_j = i$, and where indices are taken modulo $N$. The associated generalised cellular integral is defined, when it converges,  by
\begin{equation} \label{GenI}
I_{\pi}(\ba,\bb)= \idotsint\limits_{0\leq t_1 \leq \dots\leq t_{N-3}\leq 1}   f_{\pi}(\ba,\bb)\,  \omega_{\pi} \end{equation}
where we set $z_1=0$, $z_2=t_1$, \dots, $z_{N-2} = t_{N-3}$, $z_{N-1} = 1$ and $z_N=\infty$, and $\omega_{\pi}$ is the cellular integrand defined in \emph{loc.\ cit.}.  In the symmetric case, when all parameters $\ba$ and $\bb$ are equal to a non-negative integer $n$, then  one retrieves the `basic'  (or totally symmetric) cellular integrals:
\[ I_{\pi}(n) = \idotsint\limits_{0\leq t_1 \leq \dots\leq t_{N-3}\leq 1} (f_{\pi})^n \,  \omega_{\pi} \]
where  $f_{\pi}$ is defined from $f_{\pi}(\ba,\bb)$ by setting all  parameters $\ba,\bb$ equal to~$1$.

The dual configuration was defined in \emph{loc.\ cit.}\  to be the (equivalence class) of the inverse permutation 
  $\pi^{-1}$. It is  convergent if and only if $\pi$ is.    It was shown in the same paper using  Poincar\'e--Verdier duality that the family of basic cellular integrals $I_{\pi^{-1}}(n)$ are dual to the $I_{\pi}(n)$ and satisfy the dual recurrence relation. 
  It would be  interesting to prove a similar result for   generalised cellular integrals.

 Here, then,  is a definition of duality in the setting of generalised cellular integrals. 
 First of all define the  rational function dual to $f_{\pi}(\ba,\bb)$ by 
\begin{equation*}
\label{GenfDual}
f^{\vee}_{\pi}(\ba,\bb) = \prod_{i\in \Z/N\Z} \frac{ (z_i - z_{i+1})^{b_{\pi_i,\pi_{i+1}}}}{ (z_{\pi^{-1}_i} - z_{\pi^{-1}_{i+1}})^{a_{i,i+1}}}.
\end{equation*}
It automatically satisfies the same homogeneity equations. A candidate for the dual of the generalised cellular integral \eqref{GenI} is therefore
\begin{equation}
\label{def:dualIgeneralised}
I_{\pi}(\ba,\bb)^{\vee}=  \idotsint\limits_{0\leq t_1 \leq \dots \leq t_{N-3} \leq 1}   f^{\vee}_{\pi}(\ba,\bb) \, \omega_{\pi^{-1}}
\end{equation}
which, by construction, is  a family of generalised cellular integrals for $\pi^{-1}$. 
In the totally symmetric case, this reduces to the duality on basic cellular integrals.

Let us apply this to the permutation  $\pi \colon  (1,2,3,4,5,6,7,8) \mapsto (8,2,4,1,7,5,3,6)$ which represents the configuration called ${}^{\phantom\vee}_8\!\pi^\vee_8$  in~\cite{Br16}. The associated generalised cellular integral is precisely that defined in \eqref{I1} upon setting
\begin{gather*} 
a_{12} = a_1, \;\;  a_{23} = a_2, \;\;  \dots , \;\; a_{67}=a_6, \;\; a_{18} = a_7, \;\; a_{78} = a_5+a_1-a_3\\
b_{68} = a_5+a_6-a_8, \;\;
b_{17} = a_1+a_2+a_6-a_4-a_8, \;\;
b_{28} = a_1+a_7+a_8-a_3-a_6
\end{gather*}
where the remaining expressions for the $b_{ij}$ were given in~\eqref{b-param}. These equations  represent a general solution to the homogeneity conditions stated above. The inverse permutation  is $\pi^{-1} \colon (1,2,3,4,5,6,7,8) \mapsto (4,2,7,3,6,8,5,1)$. To bring the resulting family of integrals \eqref{def:dualIgeneralised} into more convenient form, we apply dihedral symmetries (more precisely,  $\sigma^2$ followed by the involution $\fri_1$) to  obtain:
\begin{equation}
\raisebox{1.8mm}{$\displaystyle\mathop{\bigintsss\!\dotsb\!\bigintsss}_{0<t_1<\dots<t_5<1}$}
\frac{\splitfrac{t_1^{-a_2+a_4-a_6+a_7+a_8}(t_2-t_1)^{a_1+a_2-a_4+a_6-a_8}(t_3-t_2)^{-a_2-a_3+a_4+a_5+a_8}}
{\times(t_4-t_3)^{a_2+a_3-a_8}(t_5-t_4)^{a_8}(1-t_5)^{a_5+a_6-a_8}\,\d t_1\dotsb\d t_5}}
{(1-t_1)^{a_7+1}(1-t_2)^{a_1-a_3+a_5+1}(t_5-t_2)^{a_6+1}(t_5-t_3)^{a_5+1}t_3^{a_4+1}t_4^{a_3+1}},
\label{I1dual}
\end{equation}
which are generalised cellular integrals for the configuration ${}_8\pi_8=(8,2,5,1,6,4,7,3)$.

The expected  duality between \eqref{I1} and \eqref{I1dual} should have many consequences including a duality between the corresponding contiguity relations.  In particular, for families of the form $I(\ba n)$ with $\ba=(a_1,\dots, a_8)$ fixed, the  coefficients of  zeta values in the two integrals should satisfy recurrence relations with respect to~$n$ which are dual to each other.  Secondly, their asymptotics should be reciprocal. Thirdly, the group of symmetries of both integrals should be related to each other. The last two properties are discussed in Section~\ref{sym-param} below.

The dual family \eqref{I1dual} admits a hypergeometric interpretation. Indeed, 
it was shown in \cite{Br16}  that a subfamily  of  ${}_8\pi_8$ is equivalent to a well-poised hypergeometric family of \cite{Zu03} (see also \cite{BR01, CFR08, Fi04, KR07} for related work and discussion) which gives linear forms in $1,\zeta(3), \zeta(5)$ but depends on fewer parameters.  In fact, as we shall now explain, this type of result may be extended to the fully general family \eqref{I1dual}. To this end, 
consider  the  very-well-poised hypergeometric series
\begin{align}
&
\tilde F_k(\bb)
=\tilde F_k(b_0;b_1,\dots,b_k)
\nonumber\\ &\;
=\frac{(b_0+1)!\,\prod_{j=1}^k b_j!}{\prod_{j=1}^k(b_0-b_j+1)!}
\cdot{}_{k+2}F_{k+1}\biggl(\begin{matrix}
b_0+2, \, \tfrac12b_0+2, \, b_1+1, \, \dots, \, b_k+1 \\[1pt]
\tfrac12b_0+1, \, b_0-b_1+2, \, \dots, \, b_0-b_k+2
\end{matrix}\biggm|(-1)^{k+1}\biggr)
\nonumber\\ &\;
=\sum_{\mu=0}^\infty(b_0+2\mu+2)
\frac{\Gamma(b_0+\mu+2)\prod_{j=1}^k\Gamma(b_j+\mu+1)}{\mu!\,\prod_{j=1}^k\Gamma(b_0-b_j+\mu+2)}
(-1)^{(k+1)\mu},
\label{eq:10}
\end{align}
where $b_0,b_1,\dots,b_k$ are non-negative integers satisfying $b_0\ge2b_i$ for $1\leq i\leq k$, and $k\ge4$.  These series are known to represent $\mathbb Q$-linear forms in $1$ and zeta values $\zeta(i)$ where $2\le i\le k-2$ and $i\equiv k\bmod2$  by \cite[Sect.~3]{Zu03}.

The  case for general $k$ will be discussed elsewhere, but here we shall only focus on the case $k=7$  which represents linear forms in $\mathbb Q+\mathbb Q\zeta(3)+\mathbb Q\zeta(5)$. 
Consider  the  arithmetic renormalisation of $\tilde F_7(\bb)$ defined by:
\begin{align}
F_7(\bb)
&=\frac{\splitfrac{(b_0-b_1-b_6)!\,(b_0-b_1-b_7)!\,(b_0-b_2-b_7)!}
{\times(b_0-b_3-b_5)!\,(b_0-b_4-b_5)!\,(b_0-b_4-b_6)!}}{b_2!\,b_3!}
\cdot\tilde F_7(\bb)
\label{eq-nm}
\\
&\in\mathbb Q+\mathbb Q\zeta(3)+\mathbb Z\zeta(5)
\nonumber
\end{align}
This series is intrinsically linked with the integral $I(\ba)$ via   results in the literature connecting very-well-poised series \eqref{eq:10} with multiple integrals (known as integrals of Sorokin type and of Vasilyev type). In more detail, 
by generating  an appropriate change of variables 
out of automorphisms of $\mathcal{M}_{0,k}$ for $4\leq k \leq 8$ akin to \cite[proof of Proposition~7.2]{Br16}, and combining with \cite[Theorem~2]{Zl02} and \cite[Theorem~5]{Zu03} we can identify the 8-parameter families $F_7(\bb)$ with  $_8\pi_8$ via
\[
F_7(\bb)
=\raisebox{0.1mm}{$\displaystyle\mathop{\bigints\!\!\!\dotsb\!\!\!\bigints}_{0<t_1<\dots<t_5<1}$}
\frac{\splitfrac{t_1^{b_0-b_2-b_7}(t_2-t_1)^{b_0-b_1-b_6}(t_3-t_2)^{b_0-b_4-b_5}(t_4-t_3)^{b_5}}
{\times(t_5-t_4)^{b_0-b_3-b_5}(1-t_5)^{b_0-b_4-b_6}\,\d t_1\dotsb\d t_5}}
{\splitfrac{(1-t_1)^{b_0-b_6-b_7+1}(1-t_2)^{b_0-b_1-b_4+1}(t_5-t_2)^{b_0-b_5-b_6+1}}
{\times(t_5-t_3)^{b_0-b_3-b_4+1}t_3^{b_3+1}t_4^{b_0-b_2-b_3+1}}}.
\]
Using the following formulae (which will follow naturally from the symmetric parametrisation considered in the following section) 
\begin{gather*}
b_0=a_2+a_3+a_4, \;\; b_1=-a_1+a_3+a_4, \;\; b_2=a_2, \;\; b_3=a_4, \;\; b_4=a_2+a_3-a_5, \\ b_5=a_2+a_3-a_8, \;\; b_6=a_4-a_6+a_8, \;\; b_7=a_2+a_3+a_6-a_7-a_8
\end{gather*}
we see that $F_7(\bb)$ coincides with~\eqref{I1dual}. This concludes the interpretation of the family of cellular integrals \eqref{I1} as being dual to very-well-poised hypergeometric series, and points to a general theory which we will not discuss further here.

\section{Symmetric parameters and their duals}
\label{sym-param}

Based on the graph-theoretic interpretation of the group $\frG$ and the table in Section~\ref{sec:Graph}, we are emboldened to introduce new parameters $s_0, s_1,\ldots, s_7$:
\begin{alignat*}{4}
a_1 &= s_1+s_2 , &\;\;  a_2 &= s_0-s_2 , &\;\;  a_3 &= s_2+s_3 , &\;\; a_4 &= s_0-s_3 ,
\\
a_5 &= s_3+s_4 , &\;\;  a_6 &= s_5+s_6 , &\;\; a_7 &= s_6+s_7 , &\;\; a_8 &= s_3+s_5 .
\end{alignat*}
We find that with this parametrization, the group $\frG$ acts on $\ba$ via the group $\frS_7$ which permutes $s_1,\dots, s_7$ and fixes $s_0$. %Furthermore, we  easily verify that the discriminant polynomial $D\in \mathbb{Q}[s_0,s_1,\ldots, s_7]$ is also invariant under this group action. 
The inverse transformation
\begin{gather*}
s_0 = \tfrac12(a_2+a_3+a_4) , \;\;  s_1 = a_1+\tfrac12(a_2-a_3-a_4) , \;\;  s_2 = \tfrac12(a_3+a_4-a_2), \\
s_3 = \tfrac12(a_2+a_3 -a_4) , \;\;  s_4 = a_5+\tfrac12(a_4-a_2-a_3) , \;\; s_5 = a_8+\tfrac12(a_4-a_2-a_3) , \\
s_6 = a_6-a_8+\tfrac12(a_2+a_3-a_4) , \;\; s_7 = a_7+a_8-a_6+\tfrac12(a_4-a_2-a_3)
\end{gather*}
shows that, in general, the symmetric parameters  $s_i$ are half-integers. 

With these new parameters, we find that the multiset $\bh$%
\footnote{Viewed as a set of  linear forms in generic parameters $a_i$ or $s_i$, $\bh$ is simply a set, but after specialising the values of the parameters to integers, $\bh$ becomes a set of integers with multiplicities.}
of 28 hyperplanes $h_{1}, \allowbreak\dots,\allowbreak h_{28}$ is precisely given by the union of two multisets
\[ 
\bh'=\{s_i+s_j:1\leq i<j \leq 7\} \quad\text{and}\quad \bh''=\{s_0-s_i:1\leq i \leq 7\}
\]
corresponding to the two orbits under the action of $\frG\cong \Sigma_7$.
Furthermore, if we combine all entries of $\bh$ into a symmetric $7\times7$ matrix
\begin{equation} 
H=(H_{ij})_{1\le i,j\le 7},
\qquad\text{where}\quad
H_{ij}=\begin{cases}
s_0-s_i &\text{if}\; i=j, \\[2.5pt]
s_i+s_j &\text{if}\; i\ne j,
\end{cases}
\label{H-matr}
\end{equation}
then this action is identified with simultaneous row-column permutations of~$H$.

Armed with the symmetric parameters $s_i$, we are now in a position to verify that the cellular integrals \eqref{I1} have isomorphic symmetry groups to their duals \eqref{I1dual}, and have reciprocal asymptotics.

It is known in general that the  very-well-poised hypergeometric series $\tilde F_k(\bb)$ in~\eqref{eq:10} is invariant under  any permutation of the parameters $b_1,\dots,b_k$. This symmetry group   can be interpreted as an action of a permutation group on $k$ elements acting by simultaneous row-column permutations of the $k\times k$ symmetric matrix
\begin{equation*}
B=(B_{ij})_{1\le i,j\le k},
\qquad\text{where}\quad
B_{ij}=\begin{cases}
b_i &\text{if}\; i=j, \\[2.5pt]
b_0-b_i-b_j &\text{if}\; i\ne j.
\end{cases}
%\label{B-matr}
\end{equation*}
In the case of interest, namely $k=7$, we may identify  this $7\times7$ matrix $B$ with \eqref{H-matr}, and thus identify the symmetry group $\frG$ for $\eqref{I1}$ with that of its dual  
by setting 
\[
\bb=(2s_0;s_0-s_1,s_0-s_2,\dots,s_0-s_7) .
\]
This change of variables  is equivalent to the equations for $b_0,\ldots, b_7$ in terms of $a_1,\ldots, a_8$ stated at the very end of Section~\ref{sec-duality}. Thus the expected duality between \eqref{I1} and \eqref{I1dual} is indeed verified on the level of symmetry groups.

Returning now to the original normalisation \eqref{eq-nm} of well-poised hypergeometric series $F_{7}(\bb)$, the above can   equivalently be expressed by  saying that 
\[
F_7(\bb)\cdot\prod_{i\in\cF}h_i!
\]
is invariant under $\frG$, for the 17-element index set $\cF$ defined in equation~\eqref{F17}.

This change in normalisation enables us to compare  the asymptotics of $\eqref{I1}$ with that of its dual.  
More precisely, it is known that there are three asymptotic quantities $\hat\lambda_1,\hat\lambda_2,\hat\lambda_3\in\mathbb R$ assigned to the behaviour of $F_7(\bb n)$
as $n\to\infty$ (see, e.g., \cite{Zu02a,Zu03}), which are related to the real numbers \eqref{3-asymp}. 
 With appropriate  labels (for example, setting $\hat\lambda_3=\lim_{n\to\infty}F_7(\bb n)^{1/n}$) they satisfy
\[
\hat\lambda_1=\lambda_1^{-1}, \quad
\hat\lambda_2=\lambda_2^{-1}, \quad
\hat\lambda_3=\lambda_3^{-1}
\]
for any \emph{fixed} choice of~$\ba$.

The above analysis gives  good reasons to expect that the  symmetries and duality between ${}^{\phantom\vee}_8\!\pi^\vee_8$ and $_8\pi_8$ (namely, of cellular integrals \eqref{I1} and \eqref{I1dual}) extend to general $k$-fold cellular integrals. This is bolstered  by an analysis of the  (known!) cases $k=4,5$,  which 
correspond to (self-dual) rational approximations to $\zeta(2)$ and $\zeta(3)$.

\section{Proof of Theorem~\ref{th:z5}}
\label{sec:proof}

Let us look at a concrete example $\ba=(8, 16, 10, 15, 12, 16, 18, 13)$ corresponding to the following matrix~\eqref{H-matr}:
\[
\begin{pmatrix}
17 & 8 & 9 & 10 & 11 & 12 & 13 \\
8 & 16 & 10 & 11 & 12 & 13 & 14 \\
9 & 10 & 15 & 12 & 13 & 14 & 15 \\
10 & 11 & 12 & 14 & 14 & 15 & 16 \\
11 & 12 & 13 & 14 & 13 & 16 & 17 \\
12 & 13 & 14 & 15 & 16 & 12 & 18 \\
13 & 14 & 15 & 16 & 17 & 18 & 11
\end{pmatrix}
\]
The ordered version of the corresponding 28-multiset $\bh$ is
\begin{align*}
\{8, 9, 10, 10, 11, 11, 11,\, & 12, 12, 12, 12, 13, 13, 13, 13,
\\
& 14, 14, 14, 14, 15, 15, 15, 16, 16, 16, 17, 17, 18\},
\end{align*}
so that $m_1=18$, $m_2=m_3=17$ and $m_4=m_5=16$.
(If the  refined estimate on the denominators were known, as discussed in   Remark~\ref{m_i}, then in view of
\begin{gather*}
\bh'=\{8, 9, 10, 10, 11, 11, 12, 12, 12, 13, 13, 13, 14, 14, 14, 15, 15, 16, 16, 17, 18\},
\\
\bh''=\{11, 12, 13, 14, 15, 16, 17\},
\end{gather*}
the only way that it could coincide with the above selection $\boldsymbol m=\{18,17,17,16,16\}$ would be if $\ell=3$ or~$4$ in the notation of that remark.
If  the true value of $\ell$ in this case were to differ from $3$ and~$4$,  there would be a small potential gain for the worthiness of our approximations to $\zeta(5)$.)
The orbit of the vector $\ba$ under the group $\frG$ is of full order  $5040$; by computation, the gain in the denominator coming from the factorial prefactors in the action of $\frG$ is
\[
\lim_{n\to\infty}\frac{\log\Phi_n}n=34.39425186\dotsc . 
\]
For the asymptotics, we find that
\begin{gather*}
\lim_{n\to\infty}\frac{\log|I(\ba n)|}n=-66.05784567\dots,
\\
\lim_{n\to\infty}\frac{\log|I'(\ba n)|}n
=\lim_{n\to\infty}\frac{\log|Q_n\zeta(5)-P_n|}n=-31.55296934\dots,
\\
\lim_{n\to\infty}\frac{\log Q_n}n=85.08768883\dotsc.
\end{gather*}
Because
\[
-31.55296934\hdots+18+17+17+16+16-34.39425186\hdots>0,
\]
we cannot make any conclusions  about the irrationality of $\zeta(5)$.
However, in the notation
\begin{gather*}
C_0=\lim_{n\to\infty}\frac{\log|Q_n\zeta(5)-P_n|}n,
\quad
C_1=\lim_{n\to\infty}\frac{\log Q_n}n>C_0,
\\
C_2=m_1+\dots+m_5-\lim_{n\to\infty}\frac{\log\Phi_n}n,
\end{gather*}
we can apply the calculation from Section~\ref{sec2} with $c_0=C_0+C_2$, $c_1=C_1+C_2$ to conclude that the worthiness of our approximations is
\[
\gamma(\ba)=\frac{C_1-C_0}{C_1+C_2}=0.86597135\dotsc.
\]
This establishes the result in Theorem~\ref{th:z5}.

\section{(In)conclusive comments}
\label{final}

There are still some mysteries remaining behind the true arithmetic of the integrals \eqref{I1}.
We have observed numerically that the exponents $\nu_p$ in \eqref{nu_p} extracted from the group $\frG$ action are not always optimal (although any losses  are insignificant).
For example, with the choice $\ba=(15, 20, 16, 14, 18, 17, 16, 20)$ (which leads to a slightly worse worthiness exponent~$\gamma=0.85163139\dots$) we obtain $\nu_p=3$ for the primes $p$ satisfying $\frac1{19} \le \{n/p\} < \frac1{18}$;
experimentally, we find  (for $n$ up to~40) that the inclusions \eqref{sharp_incl} remain valid if one replaces this exponent by~4 in the definition of $\Phi_n(\ba)$. 
This additional gain will increase the worthiness  in this case (very moderately) to $0.85665016\dots$\,.
Such arithmetic losses%
\footnote{No losses are detected for our choice $\ba=(8, 16, 10, 15, 12, 16, 18, 13)$ in Section~\ref{sec:proof}.}
do not  mean that there exists in reality a larger  group of transformations acting on the integrals $I(\ba)$;  as explained in Section~\ref{sec:Graph},  the group $\frG$ is exhaustive. 
Nevertheless, we expect the extra savings to come from different integral representations, each  possessing their own arithmetic.
Such a phenomenon has been observed for the simultaneous rational approximations to $\zeta(2)$ and $\zeta(3)$ in~\cite{Zu14}. But what are those other representations?
The machinery in \cite{Br16} allows one to generate a very large number of alternative  forms of the integral~\eqref{I1}, which are slightly more involved but which potentially have  arithmetic of their own.
For example, a change of variables brings the integral to the form
\begin{align}
I(\ba)
&=\mathop{\bigintss\!\!\dotsb\!\!\bigintss}_{[0,1]^5\;\;}
\frac{\splitfrac{x_1^{a_2}(1-x_1)^{a_1}x_2^{a_1 + a_2}(1-x_2)^{- a_2 - a_3 + a_4 - a_6 + a_7 + 2a_8}x_3^{a_1 - a_3 + a_4 + a_8}}
{\times (1-x_3)^{a_5}x_4^{a_1 - a_3 + a_4 + a_5}(1-x_4)^{a_6}x_5^{a_4}(1-x_5)^{a_1 - a_3 + a_5}}}
{\splitfrac{(1-x_2x_3x_4x_5)^{a_1 - a_3}(1-x_1x_2)^{a_1 - a_3 - a_6 + a_7 + a_8}(1-x_2x_3)^{a_8}}{\times(1-x_3x_4)^{- a_2 - a_3 + a_4 + a_5 + a_8}(1-x_4x_5)^{a_5 + a_6 - a_8}}}
\nonumber\\ &\quad\times
\frac{x_2x_3x_4\,\d x_1\,\d x_2\,\d x_3\,\d x_4\,\d x_5}{(1-x_1x_2)(1-x_2x_3)(1-x_3x_4)(1-x_4x_5)},
\label{I-b}
\end{align}
which is similar to \eqref{I-a} except for the extra factor $(1-x_2x_3x_4x_5)^{a_1 - a_3}$ in the denominator.
Clearly, such representations \eqref{I-b} are available for every $I(\frg\ba)$, where $\frg\in\frG$, and they do not degenerate  to the form \eqref{I-a} if the  exponent of $1-x_2x_3x_4x_5$ is nonzero.

The central idea of this note was to present a new approach to the elimination of parasitic terms amongst small linear forms in zeta values by viewing them as simultaneous approximations. The  price to pay was that one must pass to the subleading asymptotic.  The example we have studied here is by no means  isolated, and extends to the case of 7-fold and higher-dimensional cellular integrals. By way of example, consider the cellular integral corresponding to $(10,  2, 4, 1, 6, 3, 8, 5, 9, 7)$ on $\mathcal{M}_{0,10}$, referred to as `vanishing in the middle' in \cite[Sect.~10.2.6, Example~(4)]{Br16}. Experimentally, it underlies a motive of rank $4$ with semi-simple pieces $\mathbb{Q}(0)$, $\mathbb{Q}(-2)$, $\mathbb{Q}(-5)$, $\mathbb{Q}(-7)$. The totally symmetric version of this family is given by the integrals
\begin{align*}
I_n
&=\idotsint\limits_{0<t_1<\dots<t_7<1}
\bigg(\frac{t_1(t_2-t_1)(t_3-t_2)(t_4-t_3)(t_5-t_4)(t_6-t_5)(t_7-t_6)(1-t_7)}{(t_3-t_1)t_3t_5(t_5-t_2)(t_7-t_2)(t_7-t_4)(1-t_4)(1-t_6)}\bigg)^n
\\ &\qquad\qquad\times
\frac{\d t_1\dotsb\d t_7}{(t_3-t_1)t_3t_5(t_5-t_2)(t_7-t_2)(t_7-t_4)(1-t_4)(1-t_6)}.
\end{align*}
On the surface, they give linear forms in $1, \zeta(2), \zeta_5, \zeta_7$ where $\zeta_5,\zeta_7$ are (linear combinations of) multiple zeta values of weights 5 and 7, respectively.
Upon closer scrutiny, we may write them as a combination of two linear forms:
\[
I_n= I'_n +  I''_n \zeta(2),
\]
where $I''_n$ is a linear form in $1,\zeta(3), \zeta(5)$ and $I'_n$ is a linear form in $1, \zeta(5), \zeta(7)$ obtained by setting $\zeta(2)=0$ in $I_n$. For example, the cases $n=0,1,2$ give
\begin{align*} 
I_0&  = \frac{75}{4}\zeta(7) - 9\zeta(5)\zeta(2),
\\ 
I_1 & = \Big(61\cdot\frac{75}{4}\zeta(7) - 300\cdot3\zeta(5) - 220\Big)
- \big(61\cdot9\zeta(5) - 300\cdot2\zeta(3) + 152\big) \zeta(2),
\\
I_2 & = \Big(52921\cdot\frac{75}{4}\zeta(7) - 261153\cdot3\zeta(5) - \frac{6021219}{32}\Big)
\\ &\qquad
- \Big(52921\cdot9\zeta(5) - 261153\cdot2\zeta(3) + \frac{535857}{4}\Big)\zeta(2), 
\end{align*}
where the linear forms $I'_n$ and $I_n''$ are also small (but not as small as $I_n$ itself). Thus, a  similar analysis to the one undertaken in this note might possibly lead to a result of the form `at least one of $\zeta(5)$ and $\zeta(7)$ is irrational'.

The idea of using the subleading asymptotics of integrals to remove parasitic periods should prove  fruitful in other contexts where it may  take a more subtle form than simply setting `$\zeta(2)=0$'. 
There are several technical limitations at the moment to execute this strategy in higher weights.
For example, calculating higher weight integrals for small values of the parameters does not seem practical with current tools.
A possible way to overcome these difficulties is to create an entirely new  theoretical machinery for the arithmetic and asymptotics of the integrals in question based on the underlying algebraic geometry and contiguity relations. 
This seems to be a good challenge already for the family ${}^{\phantom\vee}_8\!\pi^\vee_8$ analysed in this note.

From a broader perspective,   this note is possibly the first systematic exploration of a cellular family of integrals which goes beyond the previously known cases for $\zeta(2)$ and $\zeta(3)$. In the process, we have uncovered several inter-related mathematical structures which include:  duality relations, the existence of large symmetry groups, and a recursive  structure  between cellular integrals  of different weights via iterated residues. 
More subtle features include the role of the symmetric parameters $s_i$, the combinatorial structures underlying the asymptotic behaviour of the integrals \eqref{I1}, and  the arithmetic of their denominators.  We expect that all these  structures point to a general theory for  cellular integrals, which   combines arithmetic, geometric, combinatorial and analytic data. Several other promising ideas for a unifying approach to irrationality proofs may  be found in~\cite{Du18, Ke20}. 
It would seem that the interplay between these different  mathematical structures, which arguably gives the subject its appeal,  lies   at the heart of the   difficulty of finding new irrationality proofs. 

\begin{acknowledgements}
This project has received funding from the European Research Council (ERC) under the European Union's Horizon 2020 research and innovation programme (grant agreement no.~724638). Both authors wish to thank the Simons Foundation for their  invitation  to the Simons Symposium on `Periods and $L$-values of Motives II',  where the idea for this  project was conceived.

We thank Ben Green for pointing out to us the vagueness in our original statement of Theorem~\ref{th:z5}, as well as Alin Bostan and Christoph Koutschan for illuminating discussions about the current state of arts for creative telescoping algorithms \cite{WZ92,Ze91}.
We also acknowledge creative feedback from the anonymous referee which caused further improvements of the text.
\end{acknowledgements}

%==================================================

\end{document}